\title{Markov chains in smooth Banach spaces\\ and Gromov hyperbolic metric spaces}

\author{Assaf Naor\thanks{Microsoft Research. Email:  anaor@microsoft.com.}\and Yuval Peres\thanks{Departments of Statistics and Mathematics,
University of California, Berkeley. Email:
peres@stat.berkeley.edu.} \and Oded Schramm\thanks{Microsoft
Research. Email:  schramm@microsoft.com.}\and Scott
Sheffield\thanks{Department of Mathematics, University of
California, Berkeley. Email:  sheff@math.berkeley.edu.}}

\documentclass[11pt]{article}
\usepackage{hyperref}
\usepackage{amssymb}
\usepackage{amsthm}
\usepackage{amsmath}
\usepackage{fullpage}
\usepackage{epsfig}
\usepackage{color,graphics}
\newcommand\remove[1]{}

\newcommand{\Lip}{\mathrm{Lip}}

\newcommand{\N}{\mathbb{N}}
\newcommand{\R}{\mathbb{R}}
\newcommand{\e}{\varepsilon}

\newcommand{\E}{\mathbb{E}}

\def\equiv{:=}

\newtheorem{theorem}{Theorem}[section]
\newtheorem{lemma}[theorem]{Lemma}
\newtheorem{prop}[theorem]{Proposition}

\newtheorem{corollary}[theorem]{Corollary}

\newtheorem{definition}[theorem]{Definition}

\begin{document}

\maketitle
\begin{abstract}
A metric space $X$ has {\em Markov type\/} $2$, if for any
reversible finite-state Markov chain $\{Z_t\}$ (with $Z_0$ chosen
according to the stationary distribution) and any map $f$ from the
state space to $X$, the distance $D_t$ from $f(Z_0)$ to $f(Z_t)$
satisfies $\E(D_t^2) \le K^2\, t\, \E(D_1^2)$ for some
$K=K(X)<\infty$.
This notion is due to K.\,Ball (1992), who showed its importance
for the Lipschitz extension problem. However until now, only
Hilbert space (and its bi-Lipschitz equivalents) were
 known to have Markov type 2.
We show that every Banach space with modulus of smoothness of
power type $2$ (in particular, $L_p$ for $p>2$)
 has Markov type $2$; this proves a conjecture of Ball. We also show that trees,
hyperbolic groups and simply connected Riemannian manifolds of
pinched negative curvature have Markov type $2$.
Our results are applied
to settle several conjectures on Lipschitz extensions
and embeddings. In particular, we answer a question
posed by Johnson and Lindenstrauss in 1982, by showing
that for $1<q<2<p<\infty$, any Lipschitz mapping from a subset of $L_p$
to $L_q$ has a  Lipschitz extension defined on all of $L_p$.
\end{abstract}

\section{Introduction}

K.\ Ball~\cite{ball} introduced the notion of Markov type of
metric spaces, defined as follows. Recall that a Markov chain
$\{Z_t\}_{t=0}^\infty$ with transition probabilities
$a_{ij}:=\Pr(Z_{t+1}=j\mid Z_t=i)$
on the state space $\{1,\ldots,n\}$
is {\em stationary\/} if $\pi_i:=\Pr(Z_t=i)$ does not depend on $t$
and it is {\em reversible\/} if $\pi_i\,a_{ij}=\pi_j\,a_{ji}$ for
every $i,j\in\{1,\ldots,n\}$.
\begin{definition}[Ball~\cite{ball}]
Given a metric space $(X,d)$ and $p\in [1,\infty)$, we say that $X$ has
{\bf Markov type} $p$ if there exists a constant $K>0$ such that for
every stationary reversible Markov chain $\{Z_t\}_{t=0}^\infty$
on $\{1,\ldots,n\}$,
every mapping $f:\{1,\ldots,n\}\to X$ and
every time $t\in \mathbb N$,
$$
\E \, d(f(Z_t),f(Z_0))^p\le K^p\,t\,\E\, d(f(Z_1),f(Z_0))^p.
$$
(Here and throughout, we omit some parentheses
and write $\E\,x^p$ for $\E(x^p)$, etc.)
The least such $K$ is called the Markov type $p$ constant of $X$,
and is denoted $M_p(X)$.
\end{definition}

Ball introduced this concept in his profound study of the
Lipschitz extension problem~\cite{ball} (see Section 2),
and the notion of Markov type has since found
applications in the theory of bi-Lipschitz
embeddings~\cite{magen-linial-naor,ramsey}. The main theorem
in~\cite{ball} states that Lipschitz functions from a subset of a
metric space $X$ having Markov type $2$ into a Banach space with
modulus of convexity of power type $2$ (see the definition
in~\eqref{def:convexity} below) extend to Lipschitz functions
defined on all of $X$. Ball showed that $M_2(L_2)=1$, yet apart
from Hilbert space and its bi-Lipschitz equivalents,
no other metric spaces were known to have Markov type $2$.
Ball asked in~\cite{ball} whether  $L_p$ for $2<p<\infty$
has Markov type $2$,
and more generally, whether
every Banach space with modulus of smoothness of power type $2$
has Markov type $2$. Recall that $X$ has  {\bf modulus of smoothness of power type} $2$ if
for all $x,y$ in the unit sphere of $X$, we have
$\|x+\tau y\|+\|x-\tau y\|-2 \le K_S\tau^2$ for some constant $K_S=K_S(X)<\infty$
and all $\tau>0$
(see also~\eqref{def:smoothness} below).
Here we answer Ball's question positively, and prove:
\begin{theorem}\label{thm:markov type2} Every normed space with
modulus of smoothness of power type $2$ has Markov type 2. Moreover,
 for every $2\le p<\infty$, we have
$ M_2(L_p)\le 4\sqrt{p-1}$.
\end{theorem}

In conjunction with Ball's extension theorem~\cite{ball}, this
implies the following non-linear version of Maurey's extension
theorem~\cite{maurey-thesis} (see below), answering a question
posed by Johnson and Lindenstrauss in~\cite{JL}.

\begin{theorem}\label{thm:extension-soft} Let $X$ be a Banach space with modulus of
smoothness of power type $2$ and let $Y$ be a Banach space with
modulus of convexity of power type $2$. Then there exists a
constant $C=C(X,Y)$ such that for every subset $A\subseteq X$ and
every Lipschitz mapping $f:A\to Y$, it is possible to extend $f$
to $\tilde f:X\to Y$ such that $\|\tilde f\|_{\Lip}\le
C\|f\|_{\Lip}$.
\end{theorem}
In particular, for  $X=L_p$ and $Y=L_q$ with $1<q< 2< p<\infty$, we can
take $C=O(\sqrt{(p-1)/(q-1)})$ (but even the finiteness of $C$ in
this case wasn't known). Additional geometric applications of
Theorem~\ref{thm:extension-soft} will be discussed later in this
introduction.

Our methods yield  Markov type $2$ for several new classes of spaces;
in particular we have the following result, which answers
positively a question posed in~\cite{naor-schechtman}.

\begin{theorem}\label{thm:tree} There exists a universal constant
$C>0$ such that for every tree $T$ with arbitrary positive edge lengths,
$$M_2(T)\le C.$$
\end{theorem}


In this theorem, as well as in
Theorem~\ref{thm:hyperbolic} and Corollary~\ref{coro:group} below,
one can take e.g.\ $C=30$. On the other hand,
we show in Section~\ref{section:discussion}
that the infinite $3$-regular tree satisfies $M_2(T)\ge \sqrt 3$.

In fact, Theorem~\ref{thm:tree} is a particular case of the
following result which holds for arbitrary Gromov hyperbolic
spaces. One of many alternative definitions for Gromov-hyperbolic
spaces is as follows (background material on this topic can be
found in the monographs~\cite{ghys,bridson,gromov-book,vaisala}).
Let $(X,d)$ be a metric space. For $x,y,r\in X$ the {\em Gromov
product} with respect to $r$ is defined as:
\begin{eqnarray}\label{eq:product}
\langle x|y\rangle_r \equiv \frac{d(x,r)+d(y,r)-d(x,y)}{2}.
\end{eqnarray}
For $\delta\ge 0$, the metric space $X$ is said to be $\delta$-hyperbolic
if for every $x,y,z,r\in X$,
\begin{eqnarray}\label{eq:def hyperbolic}
\langle x|y\rangle_r\ge \min\left\{\langle x|z\rangle_r ,\langle
y|z\rangle_r\right\}-\delta,
\end{eqnarray}
and $X$ is {\em Gromov-hyperbolic} if it is $\delta$-hyperbolic
for some $\delta<\infty$.

For trees, $\langle x|y\rangle_r$ is precisely the distance of
the least common ancestor of $x$ and $y$ from the root $r$.
It is easy to conclude (and is well known) that
every tree is $0$-hyperbolic.
Thus, the following theorem generalizes Theorem~\ref{thm:tree}

\begin{theorem}\label{thm:hyperbolic} There exists a universal
constant $C>0$ such that for every $\delta$-hyperbolic metric
space $X$, every stationary reversible Markov chain
$\{Z_t\}_{t=0}^\infty$ on
$\{1,\ldots,n\}$,
 every $f:\{1,\ldots,n\}\to X$ and every time $t\ge 1$,
$$
\E\, d(f(Z_t),f(Z_0))^2\le C^2\,t\, \E\,
d(f(Z_1),f(Z_0))^2+C^2\,\delta^2\,(\log t)^2.
$$
\end{theorem}

Hyperbolic groups are finitely generated groups on which the word
metric is $\delta$-hyperbolic for some $0\le \delta<\infty$
(see~\cite{ghys}). Since the minimal distance in such metric
spaces is $1$, we have the following corollary of
Theorem~\ref{thm:hyperbolic}:
\begin{corollary}\label{coro:group} There exists a universal
constant $C>0$ such that for every $\delta$-hyperbolic group $G$
equipped with word metric,
$$
M_2(G)\le C(1+\delta).
$$
\end{corollary}

In Section~\ref{section:hyperbolic} we also prove the following
result:
\begin{theorem}\label{thm:manifold}
Let $X$ be an $n$ dimensional, complete, simply connected Riemannian
manifold with pinched negative sectional curvature, i.e., its
sectional curvature takes values in the interval $[-R,-r]$,
where $0<r<R<\infty$. Then
$X$ has Markov type $2$ and $M_2(X)$ can be bounded by a function
of $n,r,R$.
\end{theorem}

This paper is organized as follows. In
Section~\ref{section:history} we discuss the history and
motivation of Ball's Markov type $2$ problem, as well as the
background we require from the geometry of Banach spaces. In
Section~\ref{section:applications} we briefly describe some
applications of our results to the extension problem for H\"older
functions and bi-Lipschitz embeddings.
Section~\ref{section:smooth} deals with the behavior of Markov
chains in smooth normed spaces, and contains the solution of
Ball's Markov type $2$ problem. Section~\ref{section:hyperbolic}
deals with the case of trees and hyperbolic metric spaces.
Section~\ref{section:weak} proves that the Laakso graphs have
Markov type~2. It is also proved there that doubling spaces and
planar graphs have a weak form of Markov type~2. Finally,
Section~\ref{section:discussion} contains some open problems.

\section{Linear and non-linear type and cotype; the linear and
Lipschitz extension problems}\label{section:history}

The classical Jordan-von Neumann theorem~\cite{jordan} states
that Hilbert space is characterized among Banach spaces by the
parallelogram identity $\|x+y\|^2+\|x-y\|^2=2\|x\|^2+2\|y\|^2$.
Inductively it follows that the generalized parallelogram identity
$2^{-n}\sum_{\e_i\in
\{-1,+1\}}\|\e_1x_1+\ldots+\e_nx_n\|^2=\sum_{i=1}^n\|x_i\|^2$
encodes the rich geometric structure of Hilbert space. In the
early 1970's, the work of Dubinsky-Pe\l
czy\'nsky-Rosenthal~\cite{DPR}, Hoffmann-J\o
rgensen~\cite{hoffman}, Kwapien~\cite{kwapien},
Maurey~\cite{maurey-def} and Maurey-Pisier~\cite{maurey-pisier}
has led to the notions of (Rademacher) {\em type} and {\em
cotype}, which are natural relaxations of the generalized
parallelogram identity. A Banach space $X$ is said to have type
$p>0$ if there exists a constant $T>0$ such that for every $n$ and
every $x_1,\ldots,x_n\in X$,
$$
\frac{1}{2^n}\sum_{\e_i\in
\{-1,1\}}\Bigl\|\sum_{i=1}^n\e_ix_i\Bigr\|^p\le T^p\sum_{i=1}^n
\|x_i\|^p.
$$
The least such constant $T$ is called the type $p$ constant of
$X$, and is denoted $T_p(X)$. Similarly, $X$ is said to have
cotype $q$ if there exists a constant $C>0$ such that for every $n$
and every $x_1,\ldots,x_n\in X$,
$$
\frac{1}{2^n}\sum_{\e_i\in
\{-1,1\}}\Bigl\|\sum_{i=1}^n\e_ix_i\Bigr\|^q\ge
\frac{1}{C^q}\sum_{i=1}^n \|x_i\|^q.
$$
The least such constant $C$ is called the cotype $q$ constant of
$X$, and is denoted $C_q(X)$.

These seemingly simple parameters turn out to encode a long list
of geometric and analytic properties of a normed space $X$, and in
the past three decades the theory of type and cotype has developed
into a deep and rich theory. We refer to the survey
article~\cite{maurey-survey}, the
books~\cite{lindenstrauss-tzafriri,milman-schechtman,pisier-book,woj,diestel}
and the references therein for numerous results and techniques in
this direction.

A fundamental result of Maurey~\cite{maurey-thesis} states that
any bounded {\em linear} operator from a linear subspace of a
Banach space with type $2$ into a Banach space of cotype $2$
extends to a bounded linear operator defined on the entire space.
More precisely:

\begin{theorem}[Maurey's extension theorem] Let $X, Y$ be Banach
spaces, $Z$ a linear subspace of $X$ and $T:Z\to Y$ a bounded
linear operator. Then there exists a linear operator $\tilde
T:X\to Y$ such that $\tilde T|_Z=T$ and $\|\tilde T\|\le
T_2(X)\,C_2(Y)\,\|T\|$.
\end{theorem}

The extension problem for linear operators between Banach spaces
is a linear variant of the classical {\em Lipschitz extension
problem} which asks for conditions on a pair of metric spaces
$X,Y$ implying that every Lipschitz function defined on a subset
of $X$ taking values in $Y$ can be extended to a Lipschitz
function defined on all of $X$, with only a bounded multiplicative
loss in the Lipschitz constant. Formally, let $X,Y$ be metric
spaces and denote by $e(X,Y)$ the least constant $K$ such that for
every $Z\subseteq X$ every Lipschitz function $f:Z\to Y$ can be
extended to a function $\tilde f:X\to Y$ satisfying $\|\tilde
f\|_{\Lip}\le K\|f\|_{\Lip}$ (if no such $K$ exists we set
$e(X,Y)=\infty$). Estimating $e(X,Y)$ under various geometric
conditions on $X$ and $Y$ is a classical problem which dates back to the
1930's. We refer to ~\cite{ww,benyamini-lindenstrauss,lee-naor}
 and the references therein for an
account of known results on Lipschitz extension.

The modern theory of the Lipschitz extension problem between
normed spaces starts with the work of
Marcus-Pisier~\cite{marcus-pisier} and
Johnson-Lindenstrauss~\cite{JL}. In~\cite{JL} it is asked if there
is a non-linear analog of Maurey's extension theorem. To
investigate this question it is clearly necessary to develop
non-linear variants of type and cotype.
While  there has been substantial progress in the past 20 years
on non-linear notions of type, a
satisfactory notion of non-linear cotype remains elusive.
Enflo~\cite{enflo1,enflo2,enflo3} studied the notion of
{\em roundness of metric spaces}
and subsequently in~\cite{enflo4}, generalized
roundness to a notion which is known today as {\em Enflo type}. Let $X$
be a metric space and fix $n\in \mathbb N$. An {\em $n$-dimensional
cube\/} in $X$ is a mapping $\e\mapsto x_\e$ from $\{-1,1\}^n$ to
$X$.
$X$ is said to have Enflo type $p$ with
 constant $K$ if for every $n$ and every $n$
dimensional cube $x:\{-1,1\}^n\to X$,
$$
\sum_{\e\in\{-1,1\}^n} d(x_\e,x_{-\e})^p\le K^p \sum_{\e\sim\e'}
d(x_{\e},x_{\e'})^p,
$$
where $\e\sim\e'$ if $\|\e-\e'\|_1=2$.

In the case of normed spaces $X$, Enflo type $p$ clearly implies
Rademacher type $p$ --- this follows by considering
cubes of the form $x_\e=\sum_{i=1}^n\e_i\,y_i$,
where $y_1,\ldots,y_n\in X$. A variant of Enflo
type was introduced and studied by Bourgain, Milman and Wolfson
in~\cite{BMW} (see also~\cite{pisier-varenna}).
In~\cite{naor-schechtman}, it was shown that for a wide class
of normed spaces, Rademacher type $p$ implies Enflo type $p$.
Despite the usefulness of these notions of non-linear type to
various embedding problems, they have not yielded
 extension theorems for Lipschitz functions.

A breakthrough on the Lipschitz extension problem was obtained in
the paper of Ball~\cite{ball}, where
he introduced the notion of Markov type $p$ (which, as shown
in~\cite{naor-schechtman}, implies Enflo type $p$).

To state the main result of~\cite{ball} we need to recall the
notions of uniform convexity and smoothness of normed spaces
(see~\cite{lindenstrauss-tzafriri} for a detailed account of this
theory). Let $(X,\|\cdot\|)$ be a normed space. The {\em modulus of
uniform convexity\/} of $X$ is defined for $\e\in [0,2]$ as
\begin{eqnarray}\label{def:convexity}
\delta_X(\e)=\inf\left\{ 1-\frac{\|x+y\|}{2}:\ x,y\in X,\
\|x\|=\|y\|=1,\ \|x-y\|=\e\right\}.
\end{eqnarray}
The normed space $X$ is said to be {\bf uniformly convex} if
$\delta_X(\e)>0$ for all $\e\in (0,2]$. Furthermore, $X$ is said
to have modulus of convexity of power type $q$ if there exists a
constant $c$ such that $\delta(\e)\ge c\,\e^q$ for all $\e\in
[0,2]$. It is straightforward to check that in this case $q\ge 2$.
By Proposition 7 in~\cite{BCL} (see also~\cite{figiel}), $X$ has
modulus of convexity of power type $q$ if and only if there exists
a constant $K>0$ such that for every $x,y\in X$
\begin{eqnarray}\label{eq:two point convex}
2\,\|x\|^q+\frac{2}{K^q}\,\|y\|^q\le \|x+y\|^q+\|x-y\|^q.
\end{eqnarray}
The least $K$ for which \eqref{eq:two point convex} holds is called the
$q$-convexity constant of $X$, and is denoted $K_q(X)$.

 The
modulus of uniform smoothness of $X$ is defined for $\tau>0$ as
\begin{eqnarray}\label{def:smoothness}
\rho_X(\tau)=\sup\left\{\frac{\|x+\tau y\|+\|x-\tau y\|}{2}-1:\
x,y\in X,\  \|x\|=\|y\|=1\right\}.
\end{eqnarray}
$X$ is said to be uniformly smooth if $\lim_{\tau\to 0}
\frac{\rho_X(\tau)}{\tau} =0$. Furthermore, $X$ is said to have modulus of
smoothness of power type $p$ if there exists a constant $K$ such
that $\rho_X(\tau)\le K\tau^p$ for all $\tau>0$. It is
straightforward to check that in this case necessarily $p\le 2$.
By Proposition 7 in~\cite{BCL}, $\rho_X(\tau)\le K\tau^p$ for every
$\tau>0$ if and only if there exists a constant $S>0$ such that
for every $x,y\in X$
\begin{eqnarray}\label{eq:two point}
\|x+y\|^p+\|x-y\|^p\le 2\,\|x\|^p+2\,S^p\,\|y\|^p.
\end{eqnarray}
The least $S$ for which \eqref{eq:two point} holds is called the
$p$-smoothness constant of $X$, and is denoted $S_p(X)$.

It was shown in~\cite{BCL} (see also~\cite{figiel}) that
$K_2(L_p)\le 1/\sqrt{p-1}$ for  $1<p\le 2$,
and  $S_2(L_p)\le \sqrt{p-1}$ for $2\le
p<\infty$ (the order of magnitude of
these constants was first calculated in~\cite{hanner}).

In~\cite{figiel-pisier,figiel} (see
also~\cite{lindenstrauss-tzafriri}, Theorem 1.e.16.) it is shown
that if a Banach space $X$ has modulus of convexity of power type
$q$ then $X$ also has cotype $q$. Similarly, if $X$ has modulus of
smoothness of power type $p$ then $X$ has type $p$. Observe that
$L_1$ has cotype $2$
(see~\cite{lindenstrauss-tzafriri,milman-schechtman}), but it is
clearly not uniformly convex. There also exist spaces of type 2
which are not uniformly smooth~\cite{james,pisier-xu}, but these
spaces are much harder to construct. For all the classical
reflexive spaces, the power type of smoothness and convexity
coincide with their type and cotype, respectively.

Thus, in the context of uniformly convex and uniformly smooth
spaces, one can ask the following variant of the
Johnson-Lindenstrauss question: is it true that $e(X,Y)<\infty$
whenever $X$ is a Banach space with modulus of smoothness of power
type $2$ and $Y$ is a Banach space with modulus of convexity of
power type $2$? This is precisely the problem studied by Ball
in~\cite{ball}, who proved the following theorem:

\begin{theorem}[Ball's extension theorem] Let $X$ be a metric
space, and $Y$ a Banach space with modulus of convexity of power
type $2$. Then $ e(X,Y)\le 6\,M_2(X)\,K_2(Y)$.
\end{theorem}

In particular, Ball showed that $M_2(L_2)=1$, whence $e(L_2,L_p)\le 6/\sqrt{p-1}$ for
$1<p\le 2$. Subsequently, in a difficult paper
 that is not based on the Markov type approach,
Tsar'kov~\cite{tsarkov} showed that every Banach space $X$
with modulus of smoothness of power type $2$ satisfies $e(X,L_2)<\infty$.

 Here we prove the following result:

\begin{theorem}\label{thm:markov type} Fix $1<q\le 2$ and let $X$
be a normed space with modulus of smoothness of power-type $q$.
Then
$$
M_q(X)\le \frac{8}{(2^{q+1}-4)^{1/q}}\, S_q(X).
$$
In particular, for every $2\le p<\infty$,
$$
M_2(L_p)\le 4\sqrt{p-1}.
$$
\end{theorem}

In conjunction with Ball's extension theorem, we obtain a
non-linear analog of Maurey's extension theorem:

\begin{theorem}\label{thm:non-linear maurey} For every two Banach
spaces $X,Y$,
$$
e(X,Y)\le 24\, S_2(X)\,K_2(Y).
$$
In particular, for $2\le p<\infty$ and $1<q\le 2$,
$$
e(L_p,L_q)\le 24 \sqrt{\frac{p-1}{q-1}}.
$$
\end{theorem}

\section{Some additional geometric
applications}\label{section:applications}

In~\cite{naor} the extension problem for H\"older functions was
studied. Let $X,Y$ be metric spaces. Recall that a function
$f:X\to Y$ is $\alpha$ H\"older with constant $K$ if for every
$x,y\in X$, $d(f(x),f(y))\le Kd(x,y)^\alpha$.
Following~\cite{naor} We denote by $\mathcal{A}(X,Y)$ the set of
all $\alpha>0$ such that for all $D\subseteq X$ any $\alpha$
H\"older function $f:D\to Y$ can be extended to an $\alpha$
H\"older function defined on all of $X$ with the same constant.
Analogously, $\mathcal{B}(X,Y)$ denotes the set of all $\alpha>0$
for which there exists a constant $C>0$ such that for all
$D\subseteq X$ any function $f:D\to Y$ which is $\alpha$ H\"older
with constant $K$ can be extended to a function defined on all of
$X$ which is $\alpha$ H\"older with constant $CK$. In~\cite{naor}
the following theorem was proved, which shows that the isometric
and isomorphic extension problems for H\"older functions between
$L_p$ spaces exhibit a phase transition: for $1<p,q\le 2$ we have
$\mathcal{A}(L_p,L_q)=\left(0,p/q^*\right]$ while
$\mathcal{B}(L_p,L_q)=\left(0,p/2\right]$, where $q^*=q/(q-1)$
(note that for $1<q<2$ we have $p/q^*<p/2$). Additionally, for any
$p/2<\alpha\le 1$ there is an $\alpha$ H\"older function from a
subset of $L_p$ to $L_q$ which cannot be extended to an $\alpha$
H\"older function defined on all of $L_p$.

The sets $\mathcal{A}(L_p,L_q)$ were calculated in~\cite{ww} for
all values of $p,q$. It is shown there that:
$$
\mathcal{A}(L_p,L_q)=\left\{\begin{array}{ll}
\left(0,p/q^*\right],& \mathrm{if}\
1<p,q\le 2,\\
\left(0,p/q\right], & \mathrm{if}\ 1<p\le 2\le q<\infty,\\
\left(0,p^*/q\right], & \mathrm{if}\  2\le p,q<\infty,\\
\left(0,p^*/q^*\right],& \mathrm{if}\ 2\le p<\infty\ \mathrm{and}\
1<q\le 2.\end{array}\right.
$$
As noted in~\cite{naor}, the solution of Ball's Markov type $2$
problem (i.e., Theorem~\ref{thm:markov type}) completes the
computation of the sets $\mathcal{B}(L_p,L_q)$:
$$
\mathcal{B}(L_p,L_q)=\left\{\begin{array}{ll} \left(0,p/2\right],&
\mathrm{if}\
1<p,q\le 2,\\
\left(0,p/q\right], & \mathrm{if}\ 1<p\le 2\le q<\infty,\\
\left(0,2/q\right], & \mathrm{if}\  2\le p,q<\infty,\\
\left(0,1\right],& \mathrm{if}\ 2\le p<\infty\ \mathrm{and}\
1<q\le 2.\end{array}\right.
$$

\bigskip

We end with some applications of our results to the theory of
bi-Lipschitz embeddings of finite metric spaces into normed
spaces. Given two metric spaces $X,Y$ and a one-to-one mapping
$f:X\hookrightarrow Y$, its distortion is defined as
$\mathrm{dist}(f)=\|f\|_{\Lip}\cdot\|f^{-1}\|_{\Lip}$.
In~\cite{magen-linial-naor,ramsey} it was shown that if $G$ is a
finite graph with girth $g$ (i.e., the length of the shortest
closed cycle in $G$ is $g$) and average degree $\delta>2$, then
every embedding of $G$ (with the graph metric) into Hilbert space
incurs distortion at least $\frac{\delta-2}{2\delta}\sqrt{g}$. In
fact, the same proof shows that any embedding of $G$ into a metric
space $X$ of Markov type $p$ incurs distortion at least
$\frac{\delta-2}{2M_p(X)}g^{(p-1)/p}$. We thus obtain new classes
of spaces into which large girth graphs cannot embed with low
distortion. In particular, any embedding of $G$ into $L_p$ for
$p>2$, incurs at least the distortion
$\frac{\delta-2}{8\delta}\sqrt{g/p}$.

In~\cite{ramsey} the metric Ramsey problem for the Hamming cube
$\Omega_d=\{0,1\}^d$ (with the $L_1$ metric) was studied, and it
was shown that if $ A\subseteq \Omega_d$ embeds in Hilbert space
with distortion $D$, then $|A|\le |\Omega_d|^{1-c_*/D^2}$. (Here,
$c_*$ is a universal constant. This estimate is shown
in~\cite{ramsey} to be optimal up to logarithmic terms.) The proof
of this fact is heavily based on the analysis of Markov chains on
subsets of the cube, and on the fact that Hilbert space has Markov
type $2$. The same result holds true for embeddings into any
metric space $X$ of Markov type $2$, with the constant $c$
replaced by $c_*/M_2(X)^2$.

In~\cite{naor-schechtman} it was was shown that any embedding of
the Hamming cube $\{0,1\}^d$ into $L_p$, $p>2$, incurs distortion
at least $c_p\sqrt{d}$, where $c_p$ depends only on $p$. Using
Theorem~\ref{thm:markov type} it is possible to obtain the optimal
dependence of $c_p$ on $p$, namely any embedding of $\{0,1\}^d$
into $L_p$, $p>2$, incurs distortion at least
$\frac{1}{20}\max\{1,\sqrt{d/p}\}$. This follows directly from
Theorem~\ref{thm:markov type} by considering the Markov chain
corresponding to the standard random walk on the Hamming cube. The
fact that this lower bound is optimal follows from the following
embedding. Let $\e_1,\ldots\e_d$ be i.i.d. symmetric $\pm 1$
valued Bernoulli random variables and map $x\in \{0,1\}^d $ to the
random variable $Y_x\equiv \sum_{i=1}^d x_i \e_i$. It is well
known (see e.g.~\cite{latala}) that for every integer $k\le d$,
$\left(\E \left|\sum_{i=1}^k
\e_i\right|^p\right)^{1/p}=\Theta\left(\min\{k,\sqrt{pk}\}\right)
$, implying that for $p\le d$ and $x,y\in \{0,1\}^d$, $
c\sqrt{p/d}\cdot \|x-y\|_1\le \left(E|Z_x-Z_y|^p\right)^{1/p}\le
\|x-y\|_1$, where $c$ is a universal constant. Thus the embedding
$x\mapsto Z_x$ incurs distortion $\Omega(\sqrt{d/p})$. An
application of having good bounds in terms of $p$ on the $L_p$
distortion of the Hamming cube is the following.
In~\cite{lee-mendel-naor} it was shown that if $\{0,1\}^d$ embeds
into $\ell_\infty^k$ with distortion $D$ then $k\ge
2^{\Omega(d/D^2)}$. This fact easily follows from the above
discussion via an argument similar to that of~\cite{diamond}: By
H\"older's inequality $\ell_\infty^k$ and $\ell_{p}^k$ are $O(1)$
equivalent when $p=\log k$. Thus, the fact that $\{0,1\}^d$ embeds
into $\ell_\infty^k$ with distortion $D$ implies that $D\ge
\Omega(\sqrt{d/p})=\Omega(\sqrt{d/\log k})$, which simplifies to
give the claimed lower bound on $k$.

\section{Markov chains in smooth normed
spaces}\label{section:smooth}

We begin by recalling why the real
line has Markov type $2$. Let $\{Z_t\}_{t=0}^\infty$ be a
stationary reversible Markov chain on $\{1,\ldots,n\}$,
with transition matrix $A=(a_{ij})$, and stationary distribution $\pi$.
The Markov type $2$ inequality (with constant $1$) for $\R$ simply
says that for every $x_1,\ldots,x_n\in \R$, we have
$
\sum_{i,j}\pi_i \,(A^t)_{ij}\,(x_i-x_j)^2\le t
\sum_{i,j}\pi_i \,a_{ij}\,(x_i-x_j)^2
$.
In fact, the stronger inequality
\begin{equation} \label{markovr}
\sum_{i,j}\pi_i \, (A^t)_{ij}\, (x_i-x_j)^2\le \Lambda(t)
\sum_{i,j}a_{ij}\,\pi_i\,(x_i-x_j)^2
\end{equation}
holds,
where $\Lambda(t):=\frac{1-\lambda^t}{1-\lambda}=1+\lambda+\cdots+\lambda^{t-1}$
and $\lambda$ is the second largest eigenvalue of $A$. (Recall that $-1 \le \lambda<1$.)
Indeed, since $\pi_i \, (A^t)_{ij}=\pi_j \, (A^t)_{ji}$, expanding the square in
(\ref{markovr}) shows that inequality may be rewritten
\begin{equation} \label{markovr2}
\bigl\langle (I-A^t)x,x\bigr\rangle_\pi \le \Lambda(t) 
\,\bigl\langle (I-A)x,x\bigr\rangle_\pi\,,
\end{equation}
where $x$ is the column vector $(x_1,\ldots,x_n)$,
and $\langle\cdot,\cdot\rangle_\pi$ refers to
the inner product in $L^2(\pi)$.
If $x$ is an eigenvector, then  (\ref{markovr}) is clear; the general case follows by a
spectral decomposition (note that $A$ is self-adjoint on $L^2(\pi)$).
\smallskip

The proof of the following lemma is a slight modification of the
proof of Lemma 3.1 in~\cite{ball}.
Let $L_q(X)$ denote the collection of Borel measurable $f:[0,1]\to X$
with $\E\|f\|^q=\int_0^1 \|f\|^q<\infty$. It is a Banach space with
norm $\bigl(\E\|f\|^q\bigr)^{1/q}$.
(Of course, the choice of the interval $[0,1]$ as the domain for
$f$ is rather arbitrary. It may be replaced with any probability
space equivalent to $[0,1]$.)

\begin{lemma}\label{lem:variable} Fix $1<q\le 2$ and let $Z\in L_q(X)$.
 Then
$$
\E\| Z\|^q\le \|E Z\|^q+\frac{S_q(X)^q}{2^{q-1}-1}\cdot\E \|Z-\E
Z\|^q.
$$
\end{lemma}
\begin{proof}
Let $\theta\ge 0$ be the largest constant such that for every
$Z\in L_q(X)$,
$$
\theta\left(\E\| Z\|^q-\|E Z\|^q\right)\le \E\|Z-\E Z\|^q. $$
 Our goal is to show that $\theta\ge
(2^{q-1}-1)S_q(X)^{-q}$. To this end, fix $\e>0$ and $Z\in L_q(X)$
such that
$$
(\theta+\e)\left(\E\| Z\|^q-\|E Z\|^q\right)> \E\|Z-\E Z\|^q.
$$
By the definition~\eqref{eq:two point} of $S_q(X)$ applied to the vectors
$x=(Z+\E Z)/2$ and $y=(Z-\E Z)/2$, we have the pointwise inequality:
$$
\|Z\|^q+\|\E Z\|^q\le 2\Bigl\|\frac12 Z+\frac12\E
Z\Bigr\|^q+2\,S_q(X)^q\Bigl\|\frac12Z-\frac12\E Z\Bigr\|^q.
$$
Taking expectations, we find that
\begin{eqnarray*}
\frac{1}{\theta+\e} \E \|Z-\E Z\|^q&<& \E \|Z\|^q-\|E Z\|^q\\
&\le&  2\left(\E\Bigl\|\frac12 Z+\frac12\E
Z\Bigr\|^q-\Bigl\|\E\left(\frac12 Z+\frac12\E
Z\right)\Bigr\|^q\right)+2\,S_q(X)^q\,\E\Bigl\|\frac12Z-\frac12\E
Z\Bigr\|^q\\
&\le&\frac{2}{\theta}\,\E\Bigl\|\frac12Z-\frac12\E
Z\Bigr\|^q+2\,S_q(X)^q\,\E\Bigl\|\frac12Z-\frac12\E Z\Bigr\|^q.
\end{eqnarray*}
It follows that $\frac{2^q}{\theta+\e}\le
\frac{2}{\theta}+2S_q(X)^q$. Letting $\e$ tend to zero, and
simplifying, yields the required result.
\end{proof}

The following theorem, first proved in \cite{pisier} (without the
explicit constant), is a simple corollary of
Lemma~\ref{lem:variable} (see also Proposition 3.3
in~\cite{ball}).

\begin{theorem}[Pisier]\label{thm:martingale} \, Fix $1<q\le 2$ and let
$\{M_k\}_{k=0}^n\subseteq L_q(X)$ be a martingale in $X$. Then
$$
\E \|M_n-M_0\|^q\le
\frac{S_q(X)^q}{2^{q-1}-1}\cdot\sum_{k=0}^{n-1}\E
\|M_{k+1}-M_k\|^q.
$$
\end{theorem}
\begin{proof} Assume that $\{M_k\}_{k=0}^n$ is a martingale with
respect to the filtration $\mathcal{F}_0\subseteq
\mathcal{F}_1\subseteq \cdots\subseteq \mathcal{F}_{n-1}$; 
that is,
$\E\left(M_{i+1} |\mathcal{F}_{i}\right)=M_i$
for $i=0,1,\dots,n-1$. By
Lemma~\ref{lem:variable} with conditioned expectation replacing expectation,
$$
\E\left(\left.\|M_n-M_0\|^q\right|\mathcal{F}_{n-1}\right)\le
\left\|M_{n-1}-M_0\right\|^q+\frac{S_q(X)^q}{2^{q-1}-1}\cdot\E(\|M_{n}-M_{n-1}\|^q|\mathcal{F}_{n-1}).
$$
Hence
$$
\E\|M_n-M_0\|^q\le
\E\left\|M_{n-1}-M_0\right\|^q+\frac{S_q(X)^q}{2^{q-1}-1}\cdot\E\|M_{n}-M_{n-1}\|^q,
$$
and the required inequality follows by induction.
\end{proof}

The following lemma is motivated by the continuous martingale
decompositions of Stochastic integrals constructed in~\cite{lyons}.

\begin{lemma}\label{lem:decomposition} Let $X$ be a normed space, $\{Z_t\}_{t=0}^\infty$ a stationary reversible Markov chain
on $\{1,\ldots,n\}$ and $f:\{1,\ldots,n\}\to X$. Then for every
$t\in \mathbb N$ there are two $X$-valued martingales
$\{M_s\}_{s=0}^t$ and $\{N_s\}_{s=0}^{t}$ (with respect to two
different filtrations) with the following properties:
\begin{enumerate}
\item For every $1\le s\le t-1$ we have that
\begin{eqnarray}\label{eq:identity}
f(Z_{s+1})-f(Z_{s-1})=(M_{s+1}-M_s)-(N_{t-s+1}-N_{t-s}).
\end{eqnarray}
\item For every $0\le s\le t-1$ and $q\ge 1$,
\begin{eqnarray}\label{eq:bound}
\max\left\{\E\|M_{s+1}-M_s\|^q,\E\|N_{s+1}-N_s\|^q\right\}\le
2^q\,\E\|f(Z_1)-f(Z_0)\|^q.
\end{eqnarray}
\end{enumerate}
\end{lemma}

\begin{proof}
Let $A=(a_{ij})$ be the transition matrix of $Z_t$,
and let $\pi_i:=\Pr(Z_0=i)$.
 Define
$$
Lf(i)=\sum_{j=1}^n a_{ij}[f(j)-f(i)]=\sum_{j=1}^n a_{ij}f(j)-f(i).
$$
Then
$$
\E (f(Z_s)|Z_0,\ldots,Z_{s-1})=
\E (f(Z_s)|Z_{s-1})=
Lf(Z_{s-1})+Z_{s-1}.
$$
Since $\{Z_s\}_{s=0}^\infty$ is reversible, we also have that for
every $0\le s<t$,
$$
\E(f(Z_s)|Z_{s+1},\ldots,Z_{t})=Lf(Z_{s+1})+f(Z_{s+1}).
$$

It follows that if we define $M_0=f(Z_0)$ and for $s\ge 1$:
$$
M_s :=f(Z_s)-\sum_{r=0}^{s-1}Lf(Z_r).
$$
then
$$
\E(M_s|Z_0,\ldots,Z_{s-1})=Lf(Z_{s-1})+f(Z_{s-1})-\sum_{r=0}^{s-1}Lf(Z_r)=M_{s-1},
$$
i.e., $\{M_s\}_{s=0}^\infty$ is a martingale with respect to the
natural filtration induced by $\{Z_s\}_{s=0}^\infty$.

Now, define 
$N_0:=f(Z_t)$ and for $1\le s\le t$
$$
N_s:=f(Z_{t-s})-\sum_{r=t-s+1}^tLf(Z_r).
$$
Then for $s\ge 1$,
$$
\E (N_s|Z_{t-s+1},\ldots,Z_t)=N_{s-1},
$$
i.e., $\{N_s\}_{s=0}^t$ is a martingale with respect to the
natural filtration induced by $Z_t,Z_{t-1},\ldots,Z_0$ (in
probabilistic terminology, $\{N_{t-s}\}_{s=0}^t$ is a {\em reverse
martingale}).

The identities
\begin{equation*} 
M_{s+1}-M_s=f(Z_{s+1})-f(Z_s)-Lf(Z_s),
\qquad\quad
N_{s+1}-N_{s}=f(Z_{t-s-1})-f(Z_{t-s})-Lf(Z_{t-s}).
\end{equation*}
imply~\eqref{eq:identity}.
To prove \eqref{eq:bound} observe that
for every $s\ge 0$ and $q\ge 1$,
\begin{eqnarray*}
\E
\|Lf(Z_s)\|^q=\sum_{i=1}^n\pi_i\Bigl\|\sum_{j=1}^na_{ij}[f(j)-f(i)]\Bigr\|^q\le
\sum_{i=1}^n\sum_{j=1}^n\pi_ia_{ij}\|f(j)-f(i)\|^q=\E
\|f(Z_1)-f(Z_0)\|^q.
\end{eqnarray*}
Therefore, 
\begin{multline*}
\E\|M_{s+1}-M_s\|^q
=
\E\|f(Z_{s+1})-f(Z_s)-Lf(Z_s)\|^q
\\
\le
2^{q-1}\,\E\|f(Z_{s+1})-f(Z_s)\|^q+2^{q-1}\,\E\|Lf(Z_s)\|^q\le
2^q\,\E\|f(Z_1)-f(Z_0)\|^q,
\end{multline*}
and similarly, $\E\|N_{s+1}-N_s\|^q\le 2^q\,\E\|f(Z_1)-f(Z_0)\|^q$.
\end{proof}

\begin{proof}[Proof of Theorem~\ref{thm:markov type}] Let $\{Z_s\}_{s=0}^t$, $f$, $\{M_s\}_{s=0}^t$ and
$\{N_s\}_{s=0}^t$ be as in Lemma~\ref{lem:decomposition}. Assume
first that $t$ is even, and write $t=2m$. Summing the
identity \eqref{eq:identity} over $s=1,3,5,\ldots, 2m-1$ we get
$$
f(Z_t)-f(Z_0)=\sum_{k=1}^{t/2}(M_{2k}-M_{2k-1})-\sum_{k=1}^{t/2}(N_{2k}-N_{2k-1}).
$$
By Theorem~\ref{thm:martingale} (applied to the
martingales $\sum_{k=1}^s(M_{2k}-M_{2k-1})$
and $\sum_{k=1}^s(N_{2k}-N_{2k-1})$),
 we conclude that
\begin{eqnarray*}
\E \|f(Z_t)-f(Z_0)\|^q&\le &2^{q-1}\E\,
\Bigl\|\sum_{k=1}^{t/2}(M_{2k}-M_{2k-1})\Bigr\|^q+2^{q-1}\E\,
\Bigl\|\sum_{k=1}^{t/2}(N_{2k}-N_{2k-1})\Bigr\|^q\\
&\le&
\frac{2^{q-1}S_q(X)^q}{2^{q-1}-1}\,\sum_{k=1}^{t/2}\bigl(\E\,
\|M_{2k}-M_{2k-1}\|^q+\E\,\|N_{2k}-N_{2k-1}\|^q\bigr)\\
&\le& \frac{2^{q-1}S_q(X)^q}{2^{q-1}-1}\, \frac{t}{2}
\, 2^{q+1}\E\, \|f(Z_1)-f(Z_0)\|^q.
\end{eqnarray*}
When $t$ is odd, apply the above reasoning at time $t-1$, to get
\begin{eqnarray*}
\E\, \|f(Z_{t})-f(Z_0)\|^q&\le& 2^{q-1}\E\,
\|f(Z_{t-1})-f(Z_0)\|^q+2^{q-1}\E\, \|f(Z_{t})-f(Z_{t-1})\|^q\\
&\le&
 \left(\frac{2^{3q-2}S_q(X)^q}{2^{q-1}-1} 
(t-1)+2^{q-1}\right)\E\, \|f(Z_1)-f(Z_0)\|^q\\
&\le& \frac{2^{3q-2}S_q(X)^q}{2^{q-1}-1} \, t\, \E\,
\|f(Z_1)-f(Z_0)\|^q.
\end{eqnarray*}
\end{proof}

\medskip

There are various natural variants of the notion of Markov type.
The following theorem deals with other moment inequalities for
Markov chains in uniformly smooth normed spaces:

\begin{theorem}\label{thm:more moments} For every $p\in
(1,\infty)$ and $q\in (1,2]$ there is a constant $C(p,q)\in
(0,\infty)$ with the following properties. Let $X$ be a normed
space with modulus of smoothness of power type $q$. Then for every
reversible Markov chain on $\{1,\ldots,n\}$,
$\{Z_t\}_{t=0}^\infty$, and every $f:\{1,\ldots, n\}\to X$, if
$p\le q$, then
$$
\E\, \|f(Z_t)-f(Z_0)\|^p\le
C(p,q)\,S_q(X)^q\,t\,\E\,\|f(Z_1)-f(Z_0)\|^p,
$$
i.e., $X$ has Markov type $p$. If $p>q$, then
$$
\E\, \|f(Z_t)-f(Z_0)\|^p\le
C(p,q)\,S_q(X)^p\,t^{p/q}\,\E\,\|f(Z_1)-f(Z_0)\|^p.
$$
\end{theorem}

\begin{proof} Observe that by the definition of $S_q(X)$,
for every $\tau>0$, $\rho_X(\tau)\le
S_q(X)^q\,\tau^q$.
Since $\rho_X(\tau)\le\tau$, assuming $p<q$
and $\rho_X(\tau)\le K\tau^q$, where $K\ge 1$, we have
$\rho_X(\tau)\le K\tau^p$.
By (the proof of) Proposition 7
in~\cite{BCL}, it follows that for $1<p\le q$, $S_p(X)\le C\,
S_q(X)^{q/p}$, where $C$ depends only on $p,q$. The first
result now follows from Theorem~\ref{thm:markov type}.

Assume now that $p>q$. To prove the second assertion, observe that
by a theorem of Figiel~\cite{figiel} (combined
with~\cite[Proposition 7]{BCL}), $S_q(L_p(X))\le CS_q(X)$, where
$C$ depends only on $p,q$. Recall that a sequence of elements
$x_1,x_2,\ldots$ in a Banach space $Y$ is called a {\em monotone basic
sequence\/} if for every $a_1,a_2,\ldots\in \R$ and every integer
$n$,
$$
\Bigl\|\sum_{i=1}^n a_i x_i\Bigr\|\le \Bigl\|\sum_{i=1}^{n+1} a_i
x_i\Bigr\|.
$$
By a result of Lindenstrauss~\cite{lindenstrauss-smoothness} (see
also~\cite[Proposition 2.2]{pisier}), for a Banach space $Y$,
if $\rho_Y(\tau)\le K\tau^q$ for all $\tau$, then for every
monotone basic sequence $\{x_i\}_{i\ge 0}$ in $Y$,
\begin{eqnarray*}\label{eq:basic}
\Bigl\|\sum_{i=1}^nx_i\Bigr\|^q\le 4^q\,K\sum_{i= 1}^n \|x_i\|^q.
\end{eqnarray*}
Let $\{M_i\}_{i\ge 0}\subseteq L_p(X)$ be an $X$-valued
martingale. By convexity, it follows that $\{(M_i-M_{i-1})\}_{i\ge
1}$ is a monotone basic sequence in $L_p(X)$. Moreover,
$\rho_Y(\tau)\le S_q(Y)^q\,\tau^q$ for every normed space $Y$.
Hence, for every $n$,
$$
 \|M_n-M_0\|_{L_p(X)}^q\le
 [4CS_q(X)]^q\sum_{i=1}^{n}\|M_i-M_{i-1}\|_{L_p(X)}^q.
$$
In other words,
$$
\E\, \|M_n-M_0\|^p\le [4CS_q(X)]^p\left(\sum_{i=1}^{n}
(\E\,\|M_i-M_{i-1}\|^p)^{q/p}\right)^{p/q}
=
[4CS_q(X)]^p\,n^{p/q} \,\E\,\|M_1-M_0\|^p .
$$
We now conclude the proof exactly as in the proof of
Theorem~\ref{thm:markov type}.
\end{proof}

We now state the case $X=\R$, $q=2$ of the second part
of Theorem~\ref{thm:more moments} with explicit constants.

\begin{theorem}\label{thm:real}
Let $p\in(2,\infty)$. for every reversible stationary finite
Markov chain  $\{Z_t\}_{t=0}^\infty$ on
$\{1,\ldots,n\}$ and every $f:\{1,\ldots,
n\}\to \R$, we have for all $t\ge 1$:
$$
\left(\E |f(Z_t)-f(Z_0)|^p\right)^{1/p}\le
16\,\sqrt{p\,t}\,\left(\E|f(Z_1)-f(Z_0)|^p\right)^{1/p}.
$$
\end{theorem}

The proof follows from the above argument, using
$S_2(L_p)\le \sqrt{p-1}$ for $p>2$ (\cite{BCL}, see also~\cite{figiel}).
By considering the standard random walk on
$\{1,2,\dots,n\}$ it follows that the dependence on $p$ in
Theorem~\ref{thm:real} is optimal, up to a universal
multiplicative factor.

\section{Markov chains in trees and hyperbolic metric
spaces}\label{section:hyperbolic}

In this section we prove Theorem~\ref{thm:hyperbolic} and
Theorem~\ref{thm:manifold}. We do not attempt to optimize the
constants. In particular, in the case of trees a more careful
analysis shows that one may take $C=8$ in Theorem~\ref{thm:tree}.
Since we do not believe that this is the optimal constant, we use
rougher estimates. In Section~\ref{section:discussion} we show
that the infinite $3$-regular tree satisfies $M_2(T)\ge \sqrt 3$.

\begin{lemma}\label{lem:doob} Let $\{Z_t\}_{t=0}^\infty$ be a reversible Markov chain on
$\{1,\ldots,n\}$ and $f:\{1,\ldots,n\}\to \R$. Then for every time
$t>0$,
$$
\E\max _{0\le s\le t}[f(Z_s)-f(Z_0)]^2\le 100\,t\,
\E[f(Z_1)-f(Z_0)]^2.
$$
\end{lemma}

Naturally, the proof relies on Doob's $L^2$ maximal inequality
for submartingales
\begin{equation} \label{doo}
\E\max_{0\le s\le t} M_s^2 \le 4\, \E|M_t|^2\,.
\end{equation}
See, e.g., \cite[\S4.4]{durrett}.
\begin{proof} Let $\{M_s\}_{s=0}^t$ and
$\{N_s\}_{s=0}^t$ be as in Lemma~\ref{lem:decomposition}. Observe
that
\begin{eqnarray}
\max _{0\le s\le t}[f(Z_s)-f(Z_0)]^2
&\le&
2\, \max _{\substack{0\le s\le t\nonumber\\s\
\mathrm{even}}}[f(Z_s)-f(Z_0)]^2+2\,\max _{\substack{1\le s\le
t\\s\ \mathrm{odd}}}[f(Z_s)-f(Z_{s-1})]^2.\nonumber
\\
&\le& 2\max _{\substack{0\le s\le t\nonumber\\s\
\mathrm{even}}}[f(Z_s)-f(Z_0)]^2+
2\,\sum _{\substack{1\le s\le
t\\s\ \mathrm{odd}}}[f(Z_s)-f(Z_{s-1})]^2.\nonumber
\end{eqnarray}
Therefore,
\begin{eqnarray}\label{eq:stupid}
\E\max _{0\le s\le t}[f(Z_s)-f(Z_0)]^2
&\le& 2\,\E \max _{\substack{0\le s\le t 
\\s\
\mathrm{even}}}[f(Z_s)-f(Z_0)]^2+
(t+1)\,\E[f(Z_1)-f(Z_0)]^2.
\end{eqnarray}

For even $s$ we have as in the proof of Theorem~\ref{thm:markov
type} that
$$
f(Z_s)-f(Z_0)=\sum_{k=1}^{s/2}(M_{2k}-M_{2k-1})-\sum_{k=1}^{s/2}(N_{2k}-N_{2k-1}).
$$
Thus, Doob's inequality gives
\begin{eqnarray}\label{eq:doob}
\E\max _{\substack{0\le s\le t\\ s\
\mathrm{even}}}[f(Z_s)-f(Z_0)]^2&\le& 2\,\E \max _{0\le i\le
t/2}\left[\sum_{k=1}^{i}(M_{2k}-M_{2k-1})\right]^2+2\,\E
\max_{0\le i\le t/2}\left[\sum_{k=1}^{i}(N_{2k}-N_{2k-1})\right]^2\nonumber\\
&\le& 8\,\E\left[\sum_{k=1}^{\lfloor
t/2\rfloor}(M_{2k}-M_{2k-1})\right]^2+8\,\E\left[\sum_{k=1}^{\lfloor
t/2\rfloor }(N_{2k}-N_{2k-1})\right]^2
\nonumber\\
&=&8\sum_{k=1}^{\lfloor t/2\rfloor}\E
(M_{2k}-M_{2k-1})^2+8\sum_{k=1}^{\lfloor t/2\rfloor}\E
(N_{2k}-N_{2k-1})^2\label{line:orthogonal}\\
&\le& 32\,t\, \E[f(Z_1)-f(Z_0)]^2,\nonumber
\end{eqnarray}
where 
in~\eqref{line:orthogonal} we have used
the fact that the martingale differences are orthogonal
and for the next inequality we used~\eqref{eq:bound}.

Together with~\eqref{eq:stupid}, this concludes the proof.
\end{proof}

In what follows we use the notation of the Gromov product~\eqref{eq:product}.

\begin{lemma}\label{lem:induction}
Let $X$ be a $\delta$-hyperbolic metric space. Then
for every $m\ge 1$ and $r,x_0,\ldots,x_m\in X$
\begin{eqnarray*}
d(x_0,x_m)^2&\le& 4\max_{0\le j< m}[d(x_0,r)-d(x_j,r)]^2+
4\max_{0< j\le m}[d(x_m,r)-d(x_{j},r)]^2 +{}\\&\phantom{\le}&
{}+4\sum_{j=0}^{m-1}d(x_j,x_{j+1})^2+
16\,\delta^2\lceil\log_2 m\rceil^2.
\end{eqnarray*}
\end{lemma}

\begin{proof}
Suppose first that $m$ is  a power of $2$, $m=2^k$.
Then~\eqref{eq:def hyperbolic}
gives
$$
\langle x_0|x_m\rangle_r\ge \min\left\{\langle x_0|x_{m/2}\rangle_r ,\langle
x_{m/2}|x_m\rangle_r\right\}-\delta.
$$
Hence, induction gives,
\begin{equation}\label{eq:lr}
\langle x_0|x_m\rangle_r\ge
\min\bigl\{\langle x_i|x_{i+1}\rangle_r:i=0,1,\dots,m-1\bigr\}
- k\,\delta\,.
\end{equation}
This also holds when $m$ is not a power of two, provided we take
$k:=\lceil \log_2m\rceil$ (for instance, define $x_i:=x_m$ for $m<i\le 2^k$).
Let $j$ be the index $i\in\{ 0,1,\dots,m-1\}$ giving the minimum
in~\eqref{eq:lr}.
Then
\begin{align*}
d(x_0,x_m)
&
=d(x_0,r)+d(x_m,r)-2\,\langle x_0,x_m\rangle
\\&
\le 2\,k\,\delta+d(x_0,r)+d(x_m,r)-2\langle x_j|x_{j+1}\rangle
\\&
=
2\,k\,\delta+\bigl(d(x_0,r)-d(x_j,r)\bigr)+\bigl(d(x_m,r)-d(x_{j+1},r)
\bigr)+d(x_j,x_{j+1})
\,.
\end{align*}
This implies
$$
d(x_0,x_m)^2
\le
4\,(2\,k\,\delta)^2+4\,\bigl(d(x_0,r)-d(x_j,r)\bigr)^2
+4\,\bigl(d(x_m,r)-d(x_{j+1},r)\bigr)^2
+4\,d(x_j,x_{j+1})^2\,.
$$
The lemma follows.
\end{proof}

We can now prove Theorem~\ref{thm:hyperbolic}:

\begin{proof}[Proof of Theorem~\ref{thm:hyperbolic}] Fix $r\in X$. By
Lemma~\ref{lem:induction} for every $t\ge 1$,
\begin{eqnarray*}
\E\, d(f(Z_t),f(Z_0))^2&\le& 4\,\E\,\max_{0\le j\le
t-1}[d(f(Z_0),r)-d(f(Z_j),r)]^2 +{}\\&\phantom{\le}&{}+
4\,\E\,\max_{0\le j\le
t-1}[d(f(Z_t),r)-d(f(Z_{j+1}),r)]^2+4\,t\, \E\,
d(f(Z_1),f(Z_{0}))^2+{}\\&\phantom{\le}&{}+16\,\delta^2\lceil\log_2 t\rceil^2.
\end{eqnarray*}

By Lemma~\ref{lem:doob} applied to the function $g(i):=d(f(i),r)$
we get
$$
\E\max_{0\le j\le t-1}[d(f(Z_0),r)-d(f(Z_j),r)]^2\le 100\,t\,\E\,
[d(f(Z_1),r)-d(f(Z_0),r)]^2\le 100\,t\,\E\, d(f(Z_1),f(Z_0))^2.
$$
Similarly, since the Markov chain is reversible,
$$
\E\max_{0< j\le t}[d(f(Z_t),r)-d(f(Z_{j}),r)]^2\le 100\,t\,\E\,
d(f(Z_1),f(Z_0))^2.
$$
The proof of Theorem~\ref{thm:hyperbolic} is thus complete.
\end{proof}

\subsection{A lower bound for the Markov type constant of trees} \label{sub.low}
 The following example shows that the infinite 3-regular tree $\Gamma_3$
satisfies $M_2(\Gamma_3)\ge \sqrt{3}$. Fix an integer $h$ and let $T_h$
be the complete binary tree of depth $h$ rooted at $r$. For $z\in
T_h$ denote $|z|=d(z,r)$. Consider the following transition kernel
on $T_h$. If $z$ has three neighbors, it goes to its neighbor
closer to $r$ with probability $1/2$ and to each of the other two
neighbors with probability $1/4$. Otherwise ($z$ has a single
neighbor or $z=r$), it goes to any of its neighbors with equal
probability.
This is the transition kernel for a stationary reversible Markov
chain $\{Z_t\}_{t=0}^\infty$ with stationary distribution
$\Pr(Z_0=z)=2^{-|z|}/(h+1)$. Denote $\widetilde S_t=|Z_0|-|Z_t|$
and observe that for every positive integer $n$, conditioned on
the event $\{n\le |Z_0|\le h-n\}$, the sequence $\{\widetilde
S_t\}_{t\le n}$ has the same distribution as $\{S_t\}_{t\le n}$,
the simple random walk on $\mathbb Z$ starting at $0$. Denote
$M_n\equiv \max_{t\le n} S_t$ and $\widetilde M_n\equiv \max_{t\le
n} \widetilde S_t$. By a theorem of Pitman~\cite{pitman},
$\{2M_n-S_n\}_{n\ge 0}$ has the same distribution as
$\{S_n\}_{n\ge 0}$ conditioned on $S_n> 0$ for all $n>0$ (which is
defined as the limit as $m\to\infty$ of
$\{S_n\}_{n\ge 0}$ conditioned to hit $m$ before revisiting $0$).
Let $\hat S_n$ denote this conditioned walk, which is a Markov
chain with transition probabilities $p(x,x\pm 1)=(x\pm 1)/(2x)$.
Induction easily gives $\E\hat S_n^2=3\,n$, and therefore,
$$
\E\bigl((2\widetilde M_n-\widetilde S_n)^2\bigm| n\le |Z_0|\le h-n
\bigr)= 3\,n\,.
$$
Since $\Pr(|Z_0|\in[n,h-n])= 1-\frac{2n}{h+1}$,
we get
\begin{equation}\label{eq:Mn2}
\E(2\widetilde M_n-\widetilde S_n)^2\ge
3\Bigl(1-\frac{2n}{h+1} \Bigr)n\,.
\end{equation}
Let $v$ be the (unique) vertex in $\{Z_0,\dots,Z_n\}$ closest to
the root, so that $|Z_0|-|v|=\widetilde M_n$. Note that $v$ need
not be on the geodesic connecting $Z_0$ to $Z_n$. Nevertheless, it
is unlikely to be far from this geodesic. We now make this
precise. Conditioned on $Z_0,v,\widetilde M_n$ and $\widetilde
S_n$, the vertex $Z_n$ is distributed uniformly among the vertices
$u$ such that $|u|=\widetilde S_n$ that are descendents of $v$
(that is, the simple path in $T_h$ from $u$ to the root contains
$v$). Therefore conditioned on $Z_0,v,\widetilde M_n$ the random
variable $d(Z_0,Z_n)$ is distributed as follows:
$d(Z_0,Z_n)=2\,\widetilde M_n-\widetilde S_n -2\, k$ with
probability $2^{-k-1}$ for integers $k$ satisfying $0\le
k<\min\{\widetilde M_n,\widetilde M_n-\widetilde S_n\}$ and the
conditioned probability is $2^{-k}$ for $k=\min\{\widetilde
M_n,\widetilde M_n-\widetilde S_n\}$. Thus, setting
$\ell=\min\{\widetilde M_n,\widetilde M_n-\widetilde S_n\}$ we
have that
\begin{eqnarray*}
\E (d(Z_0,Z_n)^2 \mid Z_0,v,\widetilde M_n, \widetilde
S_n)&=&\sum_{k=0}^{\ell-1}2^{-k-1} (2\,\widetilde M_n-\widetilde
S_n -2\, k)^2+2^{-\ell}(2\,\widetilde M_n-\widetilde S_n -2\,
\ell)^2\\
&\ge& (2\,\widetilde M_n-\widetilde S_n)^2-(2\,\widetilde
M_n-\widetilde S_n)\cdot \sum_{k=0}^\infty \frac{4k}{2^k}\\&=&
(2\,\widetilde M_n-\widetilde S_n)^2-8(2\,\widetilde
M_n-\widetilde S_n).
\end{eqnarray*}
Taking expectations and applying Jensen's inequality gives
$$
\E\, d(Z_0,Z_n)^2 \ge \E (2\,\widetilde M_n-\widetilde S_t )^2-
8\, \bigl({\E(2\,\widetilde M_n-\widetilde S_n )^2}\bigr)^{1/2}\,.
$$
Now, \eqref{eq:Mn2} yields
$$
\E\, d(Z_0,Z_n)^2 \ge 3\Bigl(1-\frac{2n}{h+1} \Bigr)n -8\sqrt{3n}
\,.
$$
By considering, say, $h=n^2$, we deduce that $\sup_{h}
M_2(T_h)^2\ge 3$. It follows that $M_2(\Gamma_3)^2\ge 3$ for the
infinite $3$-regular tree $\Gamma_3$, since it contains all the
finite trees $T_h$.

\section{Embeddings in products of $\R$-trees and
proof of Theorem~\ref{thm:manifold}}

In what follows, given two metric spaces $(X,d_X)$ and $(Y,d_Y)$,
the metric space $X\times Y$ is always assumed to be equipped with
the metric $d((x,y),(x',y'))=d_X(x,x')+d_Y(y,y')$.

A metric $d$ on a space  $X$ is a {\em path metric} if
for every $x,y\in X$ there is a path in $X$ from $x$ to $y$
whose length is $d(x,y)$.

An $\R$-tree is a  path metric space $(T,d)$ such that for
every two distinct points $x,y\in T$ there is a unique simple path from
$x$ to $y$ in $T$ (see~\cite{Chi,MorganShalen}).
(Some definitions appearing in the literature also require the
metric to be complete.)
Equivalently, an $\R$-tree is a $0$-hyperbolic metric space whose metric
is a path metric.

An {\em $r$-separated} set $A$ in a metric space $(X,d)$ is a
subset $A\subseteq X$ such that $d(x,x')\ge r$ for every $x\ne x'$
in $A$. An {\em $r$-net} is a maximal $r$-separated set.
If $A\subseteq X$ is an
$r$-net, then $X\subseteq \bigcup_{a\in A} B(a,r)$. Clearly, every
metric space has an $r$-net for every $r>0$.

\begin{lemma}\label{lem:glue} Fix an integer $n$ and let $Z$ be an $n$-dimensional normed space.
Let $(X,d)$ be a metric space and $D,\e>0$. Assume that every ball
of radius $\e$ in $X$ embeds bi-Lipschitzly in $Z$ with distortion
at most $D$. Then there is an integer $N=N(n,D)$,
a constant $\Delta=\Delta(n,\e,D)<\infty$ and a mapping $F:X\to Z^N$ which is
Lipschitz with constant $\Delta$ and for every $x,y\in X$
$$
d(x,y)\le \frac{\e}{8}\ \Longrightarrow\   d(x,y)\le
\|F(x)-F(y)\|.
$$
\end{lemma}

\begin{proof} Write $X_0:=X$. Having defined $X_i\subseteq X$,
let $A_{i}$ be an $\e$-net in $X_i$. Define $X_{i+1}=X_i\setminus
\bigcup_{x\in A_{i}} B(x,\e/4)$. Observe that if $x\in X_i$ then
for all $j\in \{1,\ldots,i+1\}$ there is a point $a_j\in A_j$ such
that $d(x,a_j)\le \e$. Moreover, for $j\neq j'$,
$d(a_j,a_{j'})>\e/4$. Our assumption is that for every $x\in X$
there is a function $\psi_x:B(x,\e)\to Z$ such that $\psi_x(x)=0$
and for all $y,z\in B(x,\e)$ we have $d(y,z)\le
\|\psi_x(y)-\psi_x(z)\|\le D\,d(y,z)$. Then $\|\psi_x(a_j)\|\le
\e\,D$ and the balls $B\bigl(\psi_x(a_j),\e/8\bigl) \subseteq
B(0,\e \,(D+1/8))\subseteq Z$ are disjoint. Comparison of volumes
therefore gives $i+1\le (8D+1)^n$.
 We have shown that there exists an
integer  $N\le (8D+1)^n$ and disjoint subsets $A_1,\ldots,
A_N\subseteq X$ such that for each $j$, if $a,b\in A_j$ then
$d(a,b)> \e$ and $X\subseteq \bigcup_{j=1}^N\bigcup_{a\in A_j}
B(a,\e/4)$.

For every $j\le N$ define for $a\in A_j$ and $x\in B(a,\e/2)$,
$$
f_j(x):=
\begin{cases}\psi_a(x)& \text{if $d(x,a)\le 3\e/8$},\cr
   \bigl(4-8\,\e^{-1}\,d(x,a)\bigr)\,\psi_a(x)&\text{if $3\e/8<d(x,a)\le\e/2$}.
\end{cases}
$$
%
We also set $f_j$ to be $0$ on $X\setminus\bigcup_{a\in A_j}
B(x,\e/2)$. (Observe that this defines a function on $X$ since the
balls $\{B(a,\e/2)\}_{a\in A_j}$ are disjoint.) It is
straightforward to verify that $f_j$ is Lipschitz with constant
$4\,D$.

Now define $F:X\to Z^{N}$ by $F=f_1\oplus f_2\oplus\cdots\oplus
f_N$. Then $F$ is Lipschitz with constant $4\,N\,D$. 
Moreover, fix $x,y\in X$ with $d(x,y)\le \e/8$. There exists $1\le
j\le N$ and $a\in A_j$ such that $d(x,a)\le \e/4$. Hence
$d(y,a)\le d(x,y)+d(x,a)\le 3\e/8$; so that $x,y\in B(a,3\,\e\,/8)$
and
$$
\|F(x)-F(y)\|\ge \|f_j(x)-f_j(y)\|=\|\psi_a(x)-\psi_a(y)\|\ge
d(x,y).
$$
\end{proof}

In what follows, a subset $A$ of a metric space $X$ is called
{\em $\e$-dense},
if for every $x\in X$ there is an $a\in A$ such that $d(x,a)\le \e$.

\begin{corollary}\label{coro:glue} Let $c\in(0,1]$.
Let $X,Y$ be metric spaces, $Z$
an $n$-dimensional normed space, $A\subseteq X$ an $\e$ dense
subset and $\varphi:X\to Y$ a $1$ Lipschitz mapping such that for
every $a,b\in A$, $d(\varphi(a),\varphi(b))\ge c\,d(a,b)$. Assume
that every ball of radius $\frac{64 \e}{c}$ in $X$ embeds
bi-Lipschitzly in $Z$ with distortion at most $D$. Then there exists an
integer $N=N(n,D)$ and a $K=K(n,\e,c,D)$ such that $X$ embeds bi-Lipschitzly
with distortion $K$ into $Y\times Z^N$.
\end{corollary}

\begin{proof}Let $F, N$ and $ \Delta$ be as in Lemma~\ref{lem:glue}
applied with $\e$ replaced by $64\,\e/c$. Define $g:X\to Y\times Z^N$ by
$g(x)=(\varphi(x),F(x))$. Then $g$ has Lipschitz constant bounded by
$\Delta+1$. Moreover, if $x,y\in X$ are such that $d(x,y)\le
8\,\e/c$ then $d(g(x),g(y))\ge \|F(x)-F(y)\|\ge d(x,y)$. If, on the
other hand, $d(x,y)>8\,\e/c$, then take $x',y'\in A$ satisfying
$d(x,x')\le \e$ and $d(y,y')\le \e$. It follows that
\begin{eqnarray*}
d(g(x),g(y))&\ge&
d(\varphi(x'),\varphi(y'))-d(\varphi(x),\varphi(x'))-d(\varphi(y),\varphi(y'))
\\&\ge&
c\,d(x',y')-d(x,x')-d(y,y')
\\&\ge&
c\,[d(x,y)-d(x,x')-d(y,y')]-d(x,x')-d(y,y')\\&\ge&
c\,d(x,y)-2\,\e\,(c+1)\ge \frac{c}{2}\,d(x,y).
\end{eqnarray*}
\end{proof}

The following lemma is known (it follows, for example, by specializing
the results of~\cite{khue}). However, since we
could not locate a clean reference, we will include a short proof.

\begin{lemma}\label{lem:treeext}
Let $X$ be a metric space, and let
$\phi:A\to T$ be a Lipschitz map from a subset $A\subset X$ to a
complete $\R$-tree $T$. Then $\phi$ may be extended to
a map $\tilde\phi:X\to T$ that agrees with $\phi$ on $A$
and has the same Lipschitz constant as $\phi$.
\end{lemma}

\begin{proof}
With no loss of generality, assume that the Lipschitz constant of
$\phi$ is $1$.
We will use $d_X$ and $d_T$ to denote the metrics in $X$
and $T$, respectively.
Consider first the case where $X$ contains only
one point that is not in $A$. Say,  $X=A\cup\{x_0\}$.
For any point $a\in A$, let $B_a$ be the closed ball
$B\bigl(\phi(a),d_X(x_0,a)\bigr)$ in $T$.
Then any point in $\bigcap _{a\in A} B_a$ may be chosen
as $\tilde\phi(x_0)$. Thus, we have to show that this intersection
is nonempty.
Consider $a,a'\in A$. The triangle inequality in $X$
shows that the sum of the radii of the two balls $B_a$, $B_{a'}$
is at least as large as the distance between the centers.
Thus, $B_a\cap B_{a'}$ contains a point on the unique simple
path joining $\phi(a)$ and $\phi(a')$ in $T$.

We claim that $\R$-trees satisfy the following Helly-type theorem.
If $\mathcal F$ is a nonempty finite collection of convex subsets
of $T$ and every two elements $F,F'\in \mathcal F$ intersect,
$F\cap F'\ne \emptyset$, then the whole collection intersects,
$\bigcap_{F\in\mathcal F} F\ne\emptyset$. (Here, $F$ convex means
that the metric of $T$ restricted to $F$ is a path metric.)
Indeed, suppose first that $\mathcal F=\{F_1,F_2,F_3\}$. Let
$t_1\in F_2\cap F_3$, $t_2\in F_1\cap F_3$ and $t_3\in F_1\cap
F_2$. Since there are no cycles in $T$, it follows that the three
paths, one joining $t_1$ to $t_2$, one joining $t_2$ to $t_3$ and
one joining $t_1$ to $t_3$ intersect. This intersection point will
be in $F_1\cap F_2\cap F_3$. Now, suppose that $\mathcal
F=\{F_1,F_2,\dots,F_n\}$, where $n>3$. By the previous case, we
know that any two sets in the collection $\mathcal F':= \{F_2\cap
F_1,F_3\cap F_1,\dots,F_n\cap F_1\}$ intersect. Induction then
implies that the whole collection $\mathcal F'$ intersects, which
implies that $\mathcal F$ intersects.

Since balls in $T$ are clearly convex, it follows that every
finite subcollection of $\{B_a\}_{a\in A}$ intersects. To prove
that $\bigcap _{a\in A} B_a\ne\emptyset$, we must invoke
completeness. Suppose now that $B(t,r)$ and $B(t',r')$ are two
closed balls in $T$ which intersect and none of these contains the
other. Observe that $B(t,r)\cap B(t',r')$ is also a closed ball in
$T$ whose center is the unique point in $T$ at distance
$\bigl(d_T(t,t')+r-r'\bigr)/2$ from $t$ and at distance
$\bigl(d_T(t,t')+r'-r\bigr)/2$ from $t'$, and whose radius is
$\bigl(r+r'-d_T(t,t')\bigr)/2$ (here, a single point is considered
as a ball of zero radius). Thus, the intersection
$I(A'):=\bigcap_{a\in A'} B_a$, where $A'\subset A$ is finite, is
a nonempty ball. Let $r_\infty$ be the infimum radius of any such
ball, and let $a_1,a_2,\dots$ be a sequence in $A$ such that the
radius of the ball $I(\{a_1,\dots,a_n\})$ converge to $r_\infty$.
Let $c_n$ be the center of $I(\{a_1,\dots,a_n\})$ and let $r_n$ be
the radius. If $r_\infty>0$, then the above formula for the radius
of the intersection of two balls shows that $c_n\in B_a$ for every
$a\in A$ and for every $n$ such that $r_n<2\,r_\infty$. Thus,
clearly, $\bigcap_{a\in A} B_a\ne\emptyset$ in this case. In the
case $r_\infty=0$, we have $r_n\searrow0$, and it follows that
$c_n$ is a Cauchy sequence. Invoking completeness, we conclude
that the limit $c_\infty:=\lim_n c_n$ exists. It follows that
$c_\infty\in B_a$ for every $a\in A$, because
$I(\{a,a_1,\dots,a_n\})\ne\emptyset$ for every $n$. This completes
the proof in the case $X=A\cup\{x_0\}$. The general case follows
by transfinite induction.
\end{proof}

\begin{theorem}\label{thm:glue} Fix $\delta>0$, and assume that
$X$ is a $\delta$-hyperbolic metric space whose metric is a path metric.
Assume that there exists $D,\e>0$ and $n\in \mathbb N$ such that
every ball of radius $\e$ in $X$ embeds in $\R^n$ with distortion
at most
$D$. Then there exists an integer $N$ and  a $K>0$ such that $X$
embeds bi-Lipschitzly into a product of $N$ $\R$-trees with distortion
at most $K$.
\end{theorem}

The trees in the statement of the theorem are infinite degree simplicial trees
whose edge lengths may be taken as $1$.
Though we will not need this, one can actually prove this with
bounded degree trees~\cite{benjamini-schramm} of edge length $1$.

\begin{proof}
Since balls of radius $\e$ embed in $\R^n$ with distortion at most $D$,
there is an upper bound (depending only on $n$ and $D$)
for the cardinality of any $\e/2$-separated set in $X$
whose diameter is less than $\e$. In the terminology of~\cite{bonk-schramm},
this means that $X$ has bounded growth in some scale.
 By a theorem of Bonk and Schramm~\cite{bonk-schramm}, there exists an
integer $m$ such that $X$ is quasi-isometric to a subset of the
$m$ dimensional hyperbolic space $\mathbb H^m$. This means that
there are constants $a,a',b>0$ and a mapping $f:X\to \mathbb H^m$
such that for all $x,y\in X$,
$$
a\,{d(x,y)}-b\le d(f(x),f(y))\le a'\,d(x,y)+b.
$$
(Actually, one may even take $a=a'$.)
By a theorem of Buyalo and Schroeder~\cite{buyalo}, there is an
$\R$-tree $T$ such that $\mathbb H^m$ is
quasi-isometric to a subset of $T^m$. It follows that there is
some $R>0$, an $R$-net $A$ of $X$ and a bi-Lipschitz embedding
$g:A\to T^m$.
Since the completion of an $\R$-tree is an $\R$-tree,
we may as well assume that $T$ is complete.
Lemma~\ref{lem:treeext} implies then that there is
a Lipschitz extension $\varphi :X\to T^m$ of $g$.
%
By scaling the distances in $T$ we may also assume that $\varphi$ is
$1$ Lipschitz and for every $x,y\in A$,
$d(\varphi(x),\varphi(y))\ge c\,d(x,y)$ for some constant $c>0$.

Let $A'$ be an $\e/8$ net in $X$.
Set $\Psi(r):=\sup_{a\in A'} |A'\cap B(a,r)|$.
We claim that $\Psi(r)<\infty$ for every $r>0$.
This holds true for $r=\e$, because
balls of radius $\e$ in $X$ embed in $\R^n$ with
distortion at most $D$.
Now let $a_1\in A'$, and let $x\in B(a_1,r)$.
Since $X$ has a path metric, there is a point
$x'\in X$ satisfying $d(x',a_1)=d(x',x)=d(a_1,x)/2\le r/2$.
Let $a'$ be a point in $A'$ satisfying $d(a',x')< \e/8$.
Then $a'\in B(a_1,r/2+\e/8)$ and
$x\in B(a',r/2+\e/8)$.
Thus,
$$
B(a_1,r)\subseteq \bigcup_{a'\in A'\cap B(a_1,r/2+\e/8)}
B(a',r/2+\e/8)\,.
$$
Consequently,
$\Psi(r)\le \Psi(r/2+\e/8)^2$, which now
implies that $\Psi(r)<\infty$ for every $r$.

We now show that for every $r>0$ there is an $n'=n'(r)$ and a
$D'=D'(r)$, both finite, such that every ball in $X$ of radius $r$
embeds in $\R^{n'}$ with distortion at most $D'$.
Indeed, let $x\in X$.
Set $A'_x=A'_x(r):= A'\cap B(x,r)$.
Clearly, $|A'_x|\le \Psi(r+\e/8)$.
For every $a\in A'$ let $\psi_a:B(a,\e)\to\R^n$
be a bi-Lipschitz embedding with distortion at most
$D$ satisfying $\psi_a(a)=0$, and let
\begin{equation*}
\phi_a(x):=\begin{cases}
\psi_a(x)&\text{if $d(a,x)<\e/2$,}\cr
2\,(1-\e^{-1}\,d(a,x))\,\psi_a(x)&\text{if $\e/2\le d(a,x)\le\e$,}\cr
0&\text{otherwise.}
\end{cases}
\end{equation*}
Now the required bi-Lipschitz embedding from $B(x,r)$ into
$\R^{n'}$ is given by $x\mapsto \bigl(\phi_a(x)\bigr)_{a\in A'_x}$
(with a padding of zeros to make $n'$ independent from $x$). It is
immediate that this maps satisfies the requirements. Now an
application of Corollary~\ref{coro:glue} completes the proof of
the theorem (since the real line $\R$ itself is an $\R$-tree).
\end{proof}

The following corollary contains Theorem~\ref{thm:manifold}:

\begin{corollary}\label{coro:manifold}
Let $X$ be an $n$ dimensional complete simply connected Riemannian
manifold with pinched negative sectional curvature (i.e., its
sectional curvature takes values in the compact interval
$[-R,-r]\subset (-\infty,0)$). Then
there is an integer $N=N(n,r,R)$ and $D=D(n,r,R)>0$ such that $X$
embeds bi-Lipschitzly into a product of $N$ trees with distortion
$D$. In particular, by Theorem~\ref{thm:tree}, $X$ has Markov type
$2$ and $M_2(X)$ can be bounded by a function of $n,r,R$.
\end{corollary}

\begin{proof} It is a standard fact~\cite{bridson} that $X$ is $\delta$
hyperbolic (with $\delta$ proportional to $1/r$).
The fact that there
exists $\e>0$ such that all the balls of of radius $\e$ in $X$
embed bi-Lipschitzly in $\R^n$ follows from Rauch's comparison
theorem (see~\cite{cheeger} and Chapter $8_+$
in~\cite{gromov-book}). Thus the required result follows from
Theorem~\ref{thm:glue}.
\end{proof}

\section{The Laakso graphs, doubling spaces and
weak Markov type}\label{section:weak}

\begin{figure}[ht]
\bigskip
\ \centering
\input{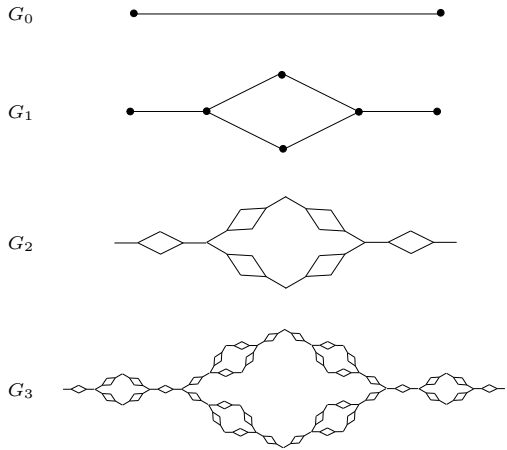}
\caption{The Laakso graphs.} \label{fig:lang}
\end{figure}

Recall that a metric space $X$ is said to be
doubling with constant $\lambda$ if every ball in $X$ can be
covered by $\lambda$ balls of half the radius. Assouad's
theorem~\cite{assouad} states that if $(X,d)$ is doubling with
constant $\lambda$ then for every $\e>0$, $(X,d^{1-\e})$ embeds in
Hilbert space with distortion $C(\lambda,\e)<\infty$. Thus $X$ has
Markov type $2-\e$ with constant depending only on $\lambda$ and
$\e$. Similarly, it is shown in~\cite{lee-mendel-naor} that if $G$ is a
planar graph equipped with the graph metric $d_G$, then for every
$\e\in (0,1)$, $(G,d_G^{1-\e})$ embeds in Hilbert space with
distortion $O(1/\sqrt{\e})$ (and the dependence on $\e$ is
optimal). Thus $G$ has Markov type $2-\e$ with constant
$O(1/\sqrt{\e})$. Is it true that a metric space with doubling
constant $\lambda$ has Markov type $2$ with constant depending
only on $\lambda$? Similarly, is it true that planar graphs have
Markov type $2$ with a universally bounded constant (this would be
a generalization of our theorem on trees)? More generally, is it
true that a Riemannian surface of genus $g$ has Markov type $2$
with constant depending only on $g$? The above embedding results
are all optimal, so a proof of the Markov type $2$ property for
such spaces cannot go through embeddings in Hilbert space. The
standard example showing that in both embedding results we {\em
must} pass to a power of the metric is the family of graphs known
as the Laakso graphs~\cite{laakso,lang}, which are planar graphs
whose graph metric is uniformly doubling, yet they do not
uniformly embed into Hilbert space. These graphs $G_k$ are defined
recursively as in Figure~\ref{fig:lang}.

\begin{prop}\label{p:lakso}
The Laakso graphs are uniformly of Markov type $2$; that is
$\sup_k M_2(G_k)<\infty$.
\end{prop}

This proposition may be viewed as some very limited indication that
doubling spaces and planar graphs have Markov type $2$.

\begin{proof} Fix $k\in\N$.
Let $r$ be the leftmost point of $G_k$, and for an arbitrary
vertex $v$ in $G_k$ set $|v|:= d(v,r)$. We claim that if
$v_0,v_1,\dots,v_t$ is a path in $G_k$, then there is some
$j\in\{0,1,\dots,t\}$ such that
\begin{equation}\label{eq:laaksopath}
d(v_0,v_t)\le \bigl||v_0|-|v_j|\bigr|+\bigl||v_j|-|v_t|\bigr|.
\end{equation}
Indeed, let $r'$ be the rightmost endpoint of $G_k$. Note that for
every vertex $v$ in $G_k$ we have $d(v,r')=|r'|-|v|$. If $v_0$ and
$v_n$ are on a geodesic path from $r$ to $r'$, then
$d(v_0,v_t)=\bigl||v_0|-|v_t|\bigr|$, and hence we may take $j=0$.
We now assume that there is no geodesic from $r$ to $r'$
containing $\{v_0,v_t\}$. Let $\beta$ be a geodesic path from
$v_0$ to $r$, and let $\gamma$ be a geodesic path from $v_t$ to
$r$. Let $u$ be the vertex in $\beta\cap\gamma$ at maximal
distance from $r$. Similarly, we define $u'$ in the same way with
$r'$ replacing $r$. It is easy to verify that  $\{u,u'\}$
separates $v_0$ from $v_t$ and separates $\{v_0,v_t\}$ from
$\{r,r'\}$.  Consequently, there is some $j\in\{0,1,\dots,t\}$
such that $v_j\in\{u,u'\}$ and that $j$
satisfies~\eqref{eq:laaksopath}.

Now suppose that $Z_0,Z_1,\dots,Z_t$ are steps of a reversible
stationary Markov chain on $\{1,\dots,n\}$ and $g$ is a mapping
from $\{0,1,\dots,n\}$ to the vertices of $G_k$.  Set
$W_j:=g(Z_j)$. For each $i=1,2,\dots,t$, consider a path in $G_k$
from $W_{i-1}$ to $W_i$ whose length is $d(W_{i-1},W_i)$. Putting
these paths together gives a path $(v_0,v_1,\dots,v_n)$ passing
through the points $W_0,\dots,W_t$ whose length $n$ is
$\sum_{i=1}^td(W_{i},W_{i-1})$. Set $j_i:= \sum_{s=1}^i
d(W_s,W_{s-1})$.  Then $v_{j_i}=W_i$. Let $j$
satisfy~\eqref{eq:laaksopath} and let $i_0\in\{1,\dots,t\}$ be
such that $j_{i_0-1}\le j\le j_{i_0}$. Then we have
\begin{eqnarray*}
d(W_0,W_t)&\le& \bigl||W_0|-|v_j|\bigr|+ \bigl||v_j|-|W_t|\bigr|\\
&\le& \bigl||W_0|-|W_{i_0-1}|\bigr|+ d(W_{i_0-1},W_{i_0})+
\bigl||W_{i_0}|-|W_t|\bigr|
\\
&\le& 2\,\bigl||W_0|-|W_{i_0-1}|\bigr|+ d(W_{i_0-1},W_{i_0})+
\bigl||W_0|-|W_t|\bigr|.
\end{eqnarray*}
Thus
\begin{eqnarray*}
d(W_0,W_t)^2&\le& 3\bigl(2\,\bigl||W_0|-|W_{i_0-1}|\bigr|\bigr)^2+
3\,d(W_{i_0-1},W_{i_0})^2+ 3\, \bigl||W_0|-|W_t|\bigr|^2
\\
&\le& 15\,\max_{0\le i\le t}\{ (|W_0|-|W_i|)^2\}+ 3\sum_{i=1}^t
d(W_{i-1},W_i)^2\,.
\end{eqnarray*}
Consequently, the proof is completed by two applications of
Lemma~\ref{lem:doob}, one with $f(Z)=|g(Z)|$ and the other with
$f(Z)=-|g(Z)|$.
\end{proof}

It turns out that in many situations a similar argument
proves a weak version of Markov type~$2$.

\begin{definition}
A metric space $X$ has weak Markov type $2$ if there is a finite
constant $C$ such that for every reversible stationary Markov
chain $Z$ on $\{1,\dots,n\}$ and every map $f:\{1,\dots,n\}\to X$,
we have
$$
\forall t\in\N,\zeta>0\qquad
\Pr(d(f(Z_0),f(Z_t))^2 \ge t\,\zeta) \le C^2\,\zeta^{-1}\,
\E(d(f(Z_0),f(Z_1))^2)\,.
$$
The least $C$ satisfying this inequality will be called the weak
Markov type $2$ constant of $X$ and will be denoted $M_2^w(X)$.
\end{definition}

Note that Chebyshev's inequality gives $M_2(X)\ge M_2^w(X)$.

\begin{theorem}\label{thm:doubling} Let $X$ be a metric space with doubling
constant $\lambda<\infty$.
Then $X$ has weak Markov type $2$, and $M^w_2(X)$
is bounded by a finite function of $\lambda$.
\end{theorem}

\begin{proof}
Let $R>0$ and let $A$ be an $R$-net in $X$.
Set $E=\{(a,a')\in A\times A:d(a,a')\le16\,R\}$.
Then $G:=(A,E)$ is a graph.
Since every ball of radius $16\,R$ in $X$ can
be covered by $\lambda^5$ balls of radius $R/2$,
it follows that the maximal degree in $G$ is
at most $\lambda^5$.
Thus, the chromatic number of $G$ is at most
$\lambda^5+1$.
Consequently, there is a partition
$A=\bigcup_{j=1}^NA_j$ of $A$ into $(16\,R)$-separated subsets
with $N\le \lambda^5+1$.
Define $f_j(x):= d(x,A_j)$, $j=1,2,\dots,N$.
Then $f_j(x)-f_j(y)\le d(x,y)$ holds for $x,y\in X$.

Now suppose that $x,y\in X$ satisfy $d(x,y)\in [3R,4R]$.
Let $a\in A$ be such that $d(a,x)\le R$ and let
$j\in\{1,\dots,N\}$ be the index such that
$a\in A_j$.
Then $f_j(x)\le d(a,x)\le R$.
Note that $d(a,y)\le d(a,x)+d(x,y)\le 5\,R$,
and since $A_j$ is $(16\,R)$-separated,
it follows that $f_j(y)=d(a,y)\ge d(y,x)-d(a,x)\ge 2\,R$.

Let $(Z_t)$ be a reversible  stationary Markov chain
on $\{1,2,\dots,n\}$ and let $f:\{1,2,\dots,n\}\to X$.
Set $W_i:=f(Z_i)$.
Fix $t\in\N$.
Let $\mathcal A$ be the event that $d(W_0,W_t)\ge 4\,R$,
let $\mathcal B_i$ be the event that $d(W_{i-1},W_{i})\ge R$,
let $\mathcal C_i$ be the event
$d(W_0,W_i)\in [3R,4R]$,
and let $\mathcal C_i^j$ be the event
$|f_j(W_0)-f_j(W_i)|\ge 2\,R$.
Clearly,
$$
\mathcal A\subseteq \left(\bigcup_{i=1}^t\mathcal
B_i\right)\bigcup \left(\bigcup_{i=1}^t\mathcal C_i\right)\,.
$$
The previous paragraph shows that
$$
\mathcal C_i\subseteq \bigcup_{j=1}^N\mathcal C_i^j\,.
$$
Consequently,
$$
\Pr(\mathcal A) \le \sum_{i=1}^t \Pr(\mathcal B_i)+
\sum_{j=1}^N\Pr\left(\max_{1\le i\le t} |f_j(W_0)-f_j(W_i)|\ge
2\,R\right)\,.
$$
The first of these sums is bounded by
$t\,\E(d(W_0,W_1)^2)\,R^{-2}$, while Lemma~\ref{lem:doob}
shows that the second sum is bounded by
$50\,N\,t\, \E(d(W_0,W_1)^2)\,R^{-2}$.
The theorem follows by choosing $R$ to satisfy $16\,R^2=t\,\zeta$.
\end{proof}

\medskip

Let $(X,d_X), (Y,d_Y)$ be metric spaces. We shall say that
$(X,d_X)$ embeds weakly into $(Y,d_Y)$ with distortion $K$ if for
every $\Delta>0$ there is a $1$ Lipschitz mapping $g_\Delta:X\to
Y$ such that if $x,y\in X$ satisfy $d_X(x,y)\ge \Delta$ then
$d_Y(g_\Delta(x),g_\Delta(y))\ge \Delta/K$. Observe that in this
case $M_2^w(X)\le K M_2^w(Y)$. Indeed fix a reversible stationary
Markov chain $Z$ on $\{1,\dots,n\}$ and a map $f:\{1,\dots,n\}\to
X$. For every time $t\in \mathbb N$ and every $\zeta>0$,
\begin{eqnarray*}
\Pr(d_X(f(Z_0),f(Z_t))^2\ge t\zeta)&\le&
\Pr\left(d_Y(g_{\sqrt{t\zeta}}(f(Z_0)),
g_{\sqrt{t\zeta}}(f(Z_t)))^2\ge \frac{t\zeta}{K^2}\right)\\&\le&
M_2^w(Y)^2K^2\zeta^{-1}\, \E \left(d_Y(g_{\sqrt{t\zeta}}(f(Z_0)),
g_{\sqrt{t\zeta}}(f(Z_t)))^2\right)\\&\le&
M_2^w(Y)^2K^2\zeta^{-1}\, \E \left(d_X(f(Z_0), f(Z_t))^2\right).
\end{eqnarray*}

It follows from the results of~\cite{lee-mendel-naor,descent} that
if $X$ is doubling with constant $\lambda$ then $X$ embeds weakly
into Hilbert space with distortion $O(\log \lambda)$. This yields
an alternative proof of Theorem~\ref{thm:doubling}, with the
concrete estimate $M_2^w(X)=O(\log \lambda)$. Moreover, the
results of~\cite{lee-mendel-naor,descent} (specifically, see the
proof of Lemma 5.2 in~\cite{lee-mendel-naor}) imply that any
planar graph embeds weakly into Hilbert space with $O(1)$
distortion. More generally, in combination with Corollary 3.15
in~\cite{lee-naor}, any Riemannian surface $\mathcal S$ of genus
$g$ embeds weakly into Hilbert space with distortion $O(g+1)$.
Thus $M_2^w(\mathcal S)=O(g+1)$.

\medskip

\section{Discussion and open problems}\label{section:discussion}
\begin{enumerate}
%
\item In~\cite{pisier-varenna} (see also a related
previous result in~\cite{BMW}) it is shown that for every $p>1$,
if a Banach space $X$ has Rademacher type $p$ then it also has
Enflo type $q$ for every $q<p$. No such result is known for Markov
type. In~\cite{naor-schechtman} it is shown that if $X$ is a UMD
Banach space (see~\cite{burkholder-survey} for details on UMD
spaces) of Rademacher type $p$, then $X$ also has Enflo type $p$.
It would be desirable to obtain a result stating that for a
certain class of Banach spaces, the notions of Rademacher type $p$
and Markov type $p$ coincide (or almost coincide). The most daring
conjecture would be that for {\em every} Banach space, Rademacher
type $p$ implies Enflo type $p$, or even Markov type $p$. This
amounts to proving that for Banach spaces of type greater than $1$
(also known as $K$-convex space. See~\cite{pisier-annals} for the
geometric and analytic ramifications of this assumption), the
Rademacher type and Enflo type (or Markov type) coincide.

One simple example of a class of spaces for which we can prove
that there is a strong connection between Rademacher type and
Markov is Banach lattices. A Banach lattice is a Banach space
$(X,\|\cdot\|)$ which is partially ordered and satisfies the
following axioms. For every $x,y,z\in X$, if $x\le y$ then $x+z\le
y+z$, and for every scalar $\alpha\in [0,\infty)$, $x\ge 0$
implies that $\alpha x\ge 0$. It is also required that for all
$x,y\in X$ there exists a least upper bound $x\vee y$ and a
greatest lower bound $x\wedge y$. For $x\in X$ denote
$|x|=x\vee(-x)$. The final requirement is that the partial
ordering is compatible with the norm in the sense that if $|x|\le
|y|$ then $\|x\|\le \|y\|$. Examples of Banach lattices are the
classical function and sequence spaces, with the point-wise
partial order. We refer to~\cite{lindenstrauss-tzafriri} for an
account of the beautiful theory of Banach lattices.

A combination of a theorem of Figiel~\cite{figiel} and a theorem
of Maurey~\cite{more-maurey} (see Theorem 1.f.1. and Proposition
1.f.17. in~\cite{lindenstrauss-tzafriri}) implies that a Banach
lattice $X$ of type $2$ can be renormed to have a modulus of
smoothness of power type $2$. Thus by Theorem~\ref{thm:markov
type} $X$ has Markov type $2$.
%
%
\item Under what conditions on a metric space does Enflo
type $p$ imply Markov type $p$?
\item
Is it true that if a metric space has
Markov type $p$ then it also has Markov type $q$ for every $q<p$?
For normed spaces this is indeed the case, by a straightforward
application of Kahane's inequality~\cite{kahane}.
%
%
\item
We conjecture that the factor of $24$ in
Theorem~\ref{thm:non-linear maurey} is redundant. In particular it
seems likely that for $2\le p<\infty$ and $1<q\le 2$,
$e(L_p,L_q)\le \sqrt{(p-1)/(q-1)}$. If true, this would be a
generalization of Kirszbraun's classical extension
theorem~\cite{kirszbraun} (see
also~\cite{ww,benyamini-lindenstrauss}).
%
%
\item
Since $L_1$ has cotype $2$ but isn't uniformly
convex, there is no known non-linear analog of Maurey's extension
theorem for $L_1$-valued mappings. In particular, it isn't known
whether $e(L_2,L_1)$ is finite or infinite.
%
%
\item What is the best Markov type $2$ constant for trees? More
precisely, define $M_2(\mathrm{tree})$ to be $\sup M_2(T)$ over
all trees $T$. (It is clear that this $\sup$ is a $\max$.) One can
show using the methods of the present paper that
$M_2(\mathrm{tree})\le 8$. The example in Subsection \ref{sub.low}
shows that
$M_2(\mathrm{tree})\ge \sqrt{3}$.
\item As discussed in Section~\ref{section:weak}, we believe that
planar graphs and doubling spaces have Markov type~2. Also, it
seems likely that CAT(0) spaces have Markov type~2
(see~\cite{bridson} for a discussion of CAT(0) spaces).

\item Say that a metric space $X$ has  {\bf maximal Markov type
$p$}, if there exists a constant $K$ such that for every finite
stationary reversible Markov chain $\{Z_t\}_{t=0}^\infty$ on
$\{1,\ldots,n\}$ and every mapping $f:\{1,\ldots,n\}\to X$, we
have
$$
\E\, \max_{1\le s\le t} d(f(Z_s),f(Z_0))^p\le K^p t\, \E\,
d(f(Z_1),f(Z_0))^p.
$$
for all $t\in \mathbb N$.
In all the cases in which we proved that a metric space $X$ has
Markov type $2$, the argument actually shows that it has {\em maximal
Markov type $2$}. This was explicit in the proofs for trees and hyperbolic spaces.
To see this in the setting of Banach
spaces with modulus of smoothness of
power type $2$, it suffices to note that Doob's $L_2$ maximal inequality
(\ref{doo}) is also valid
for a martingale $\{M_s\}_{s \ge 0}$ in a Banach space,
since by Jensen's inequality,  $\{\|M_s\|\}_{s \ge 0}$ is a submartingale for the same filtration.
We do not know whether in general, Markov type
$p$ implies maximal Markov type $p$.
\end{enumerate}

\smallskip

\noindent{\bf Acknowledgement.} We are grateful to Russ Lyons for helpful discussions
at an early stage of this work, and to Terry Lyons for sending us his paper~\cite{lyons}
 with T. S. Zhang.


\begin{thebibliography}{11}



\bibitem{assouad} P. Assouad. Plongements lipschitziens dans
$\mathbf{R}^n$. {\em Bull. Soc. Math. France} 111 429--448 (1983).

\bibitem{ball} K. Ball. Markov Chains, Riesz Transforms and
Lipschitz Maps. {\em Geom. Funct. Anal.} 2, 137--172 (1992).

\bibitem{BCL} K. Ball, E. A. Carlen and E. H. Lieb. Sharp uniform
convexity and smoothness inequalities for trace norms. {\em
Invent. math.} 115, 463--482 (1994).

\bibitem{ramsey} Y. Bartal, N. Linial, M. Mendel and A. Naor.
On metric Ramsey-type phenomena. {\em Ann. Math.}, to appear.

\bibitem{benjamini-schramm} I. Benjamini and O. Schramm.
In preparation.

\bibitem{benyamini-lindenstrauss} Y. Benyamini and J.
Lindenstrauss. Geometric Nonlinear Functional Analysis, volume 1.
Amer. Math. Soc. Coll. Publ. 48 (2000).

\bibitem{bonk-schramm} M. Bonk and O. Schramm. Embeddings of
Gromov hyperbolic spaces. {\em Geom. Funct. Anal.} 10, 266--306
(2000).

\bibitem{bridson} M. R. Bridson and A. Haefliger. Metric spaces of
non-positive curvature. Springer-Verlag, Berlin (1999).

\bibitem{BMW} J. Bourgain, V. Milman and H. Wolfson. On type of
metric spaces. {\em Trans. Amer. Math. Soc.} 294, no. 1, 295--317
(1986).




\bibitem{burkholder-survey} D. L. Burkholder. Martingales and
Singular integrals in Banach spaces. {\em Handbook of the Geometry
of Banach Spaces, volume 1}, W. B. Johnson and J. Lindenstrauss
eds., Elsevier, Amsterdam, 233--269 (2001).


\bibitem{buyalo} S. Buyalo and V. Schroeder. Embedding of
hyperbolic spaces in the product of trees. Preprint (2004),
available at http://arxiv.org/abs/math.GT/0311524.

\bibitem{cheeger} J. Cheeger and D. Ebin. Comparison Theorems in
Riemannian Geometry. North Holland, New York (1975).

\bibitem{Chi} I. M. Chiswell. Length functions and free products of groups
{\em Proc. London Math. Soc.} 42, 42--58 (1981).

\bibitem{diestel} J. Diestel, H. Jarchow and A. Tonge. Absolutely
Summing Operators. Cambridge University Press (1995).

\bibitem{DPR} E. Dubinsky, A. Pe\l czy\'nsky and H.P. Rosenthal.
On Banach spaces $X$ for which $\Pi_2(L_\infty,X)=B(L_\infty,X)$.
{\em Sudia Math.} 44, 617--648 (1972).

\bibitem{durrett} R. Durrett. Probability: Theory and Examples, second edition.
Duxbury Press, Belmont CA, (1996).

\bibitem{enflo1} P. Enflo. Topological groups in which
multiplication on one side is differentiable or linear. {\em Math.
Scand.} 24, 195--207 (1970).

\bibitem{enflo2} P. Enflo. Uniform structures and square roots in
topological groups I. {\em Israel J. Math.} 8 230--252 (1970).

\bibitem{enflo3} P. Enflo. Uniform structures and square roots in
topological groups II. {\em Israel J. Math.} 8 253--272 (1970).

\bibitem{enflo4} P. Enflo. On infinite-dimensional topological
groups. {\em S\'eminaire sur la G\'eom\'etrie des Espaces de
Banach 1977-1978}, Exp. no. 10--11, \'Ecole Polytech., Palaiseau
(1978).

\bibitem{figiel} T. Figiel. On the moduli of convexity and
smoothness. {\em Studia Math.} 56, 121--155 (1976).

\bibitem{figiel-pisier} T. Figiel and G. Pisier. S\'eries
al\'eatoires dans les espaces uniform\'ement convexes ou
uniform\'ement lisses. {\em C. R. Acad. Sci.} 279, 611--614
(1974).

\bibitem{ghys} E. Ghys and P. de la Harpe (eds.). Sur les groupes
hyperboliques d'apr\`es Mikhael Gromov. Birkh\"auser, Boston MA
(1990).

\bibitem{gromov-book} M. Gromov. Metric Structures for Riemannian
and non-Riemannian Spaces. Birkh\"auser, Boston (1999).


\bibitem{hanner} O. Hanner. On the uniform convexity of $L^p$ and
$l^p$. {\em Ark. Math.} 3, 239--244 (1956).

\bibitem{heinonen} J. Heinonen. Lectures on Analysis in Metric
Spaces. Universitext, Springer-Verlag, New York (2001).

\bibitem{hoffman} J. Hoffman-J\o rgensen. Sums of independent
Banach space valued random variables. {\em Studia Math.} 52,
159--186 (1974).

\bibitem{james} R. C. James. Nonreflexive spaces of type $2$. {\em
Israel J. Math.} 30, 1--13 (1978).

\bibitem{JL} W. B. Johnson and J. Lindenstrauss. Extensions of
Lipschitz mappings into a Hilbert space. In {\em Conference in
modern analysis and probability (New Haven, Conn., 1982)}, volume
26 0f {\em Contemp. Math.}, 189--206, Amer. Math. Soc., Providence
RI (1984).

\bibitem{gid} W. B. Johnson, J. Lindenstrauss, D. Preiss and G.
Schechtman. Lipschitz quotients from metric trees and from Banach
spaces containing $\ell_1$. {\em J. Funct. Anal.} 194, 332--346
(2002).

\bibitem{jordan} P. Jordan and J. von Neumann. On inner products in linear metric spaces. {\em Ann. of Math. (2)} 36, no. 3, 719--723
(1935).

\bibitem{kahane} J. P. Kahane. Series of Random Functions. {\em
Heath Math. Monographs}, Lexington Mass., Heath and Co. (1968).

\bibitem{khue} N. V. Khue and N. T. Nhu. Lipschitz extensions and
Lipschitz retractions in metric spaces. {\em Colloquium Math.} 45,
245--250 (1981).

\bibitem{kirszbraun} M. D. Kirszbraun. \"Uber die
zusammenziehenden und Lipschitzchen Transformationen. {\em Fund.
Math.} 22, 77--108 (1934).

\bibitem{descent} R. Krauthgamer, J. R. Lee, M. Mendel and A.
Naor. Measured descent: A new embedding method for finite metrics.
Preprint (2004).

\bibitem{kwapien} S. Kwapien. Isomorphic characterizations of
inner product spaces by orthogonal series with vector valued
coefficients. {\em Studia Math.} 44, 583--595 (1972).

\bibitem{laakso} T. J. Laakso. Ahlfors $Q$-regular spaces with arbitrary $Q>1$ admitting weak Poincar\'e inequality.
{\em Geom. Funct. Anal.} 10, no. 1, 111--123 (2000).

\bibitem{lang} U. Lang and C. Plaut. Bilipschitz embeddings of
metric spaces into space forms. {\em Geom. Dedicata}, 87 (1--3),
285--307 (2001).

\bibitem{latala} R. Lata\l a. Estimation of moments of sums of
independent real random variables. {\em Ann. Probab.} 25, no. 3,
1502--1513 (1997).

\bibitem{diamond} J. R. Lee and A. Naor. Embedding the diamond
graph in $L_p$ and dimension reduction in $L_1$. {\em Geom. Funct.
Anal.} 14, no. 4, 745--747 (2004).

\bibitem{lee-naor} J. R. Lee and A. Naor. Extending Lipschitz
functions via random metric partitions. {\em Invent. math.}, to
appear.

\bibitem{lee-mendel-naor} J. R. Lee, M. Mendel and A. Naor. Metric
structures in $L_1$: Dimension, snowflakes, and average
distortion. {\em European Journal of Combinatorics}, to appear.

\bibitem{lindenstrauss-smoothness} J. Lindenstrauss. On the
modulus of smoothness and divergent series in Banach spaces. {\em
Michigan Math. J.}, 10, 241--252 (1963).

\bibitem{lindenstrauss-tzafriri} J. Lindenstrauss and L. Tzafriri.
Classical Banach Spaces II. Springer-Verlag (1979).

\bibitem{lyons} T. J. Lyons and T. S. Zhang.
Decomposition of Dirichlet processes and its applications. {\em
Ann. Probab.} 22, 494--524 (1994).

\bibitem{magen-linial-naor} A. Magen, N. Linial and A. Naor. Girth
and Euclidean distortion. {\em Geom. Funct. Anal.} 12, no. 2,
380--394 (2002).

\bibitem{marcus-pisier} M. B. Marcus and G. Pisier.
Characterizations of almost surely continuous $p$-stable random
Fourier series and strongly stationary processes. {\em Acta Math.}
152 (3--4), 245--301 (1984).

\bibitem{maurey-def} B. Maurey. Espaces de cotype $p$, $0<p\le 2$.
{\em S\'eminaire Maurey-Schwartz 1972/73}, \'Ecole Polytech.,
Palaiseau (1973).

\bibitem{maurey-thesis} B. Maurey. Th\'eor\`ems de factorisation
pour les operateurs lin\'eaires a valuers dans un espace $L^p$.
{\em Asterisque}, 11, Soc. Math. France (1974).

\bibitem{more-maurey} B. Maurey. Type et cotype dans les espaces
munis de structures locales inconditionnelles. {\em S\'eminaire
Maurey-Schwartz 1973/74}, exposes 24--25, \'Ecole Polytech.,
Palaiseau (1974).


\bibitem{maurey-survey} B. Maurey. Type, Cotype and $K$-Convexity.
{\em Handbook of the Geometry of Banach Spaces, volume 2}, W. B.
Johnson and J. Lindenstrauss eds., Elsevier, Amsterdam,
1299--1332. (2003).

\bibitem{maurey-pisier} B. Maurey and G. Pisier. S\'eries de
variables al\'eatoire vectorielles ind\'ependants et
propri\'et\'es g\'eom\'etriques des espaces de Banach. {\em Studia
Math.} 58, 45--90 (1976).

\bibitem{milman-schechtman} V. D. Milman and G. Schechtman.
Asymptotic Theory of Finite Dimensional Normed Spaces. With an
appendix by M. Gromov. {\em Lect. Notes Math.} 1200,
Springer-Verlag, Berlin (1986).

\bibitem{MorganShalen}
J. W. Morgan, P. B. Shalen. Valuations, trees, and degenerations
of hyperbolic structures. {\em Ann. of Math. (2)} 120, 401--476
(1984).

\bibitem{naor} A. Naor. A phase transition phenomenon between the
isometric and isomorphic extension problems for H\"older functions
between $L_p$ spaces. {\em Mathematika} 48, 253--271 (2001).

\bibitem{naor-schechtman} A. Naor and G. Schechtman. Remarks on
non linear type and Pisier's inequality. {\em J. Reine Angew.
Math.} 552, 213--236 (2002).

\bibitem{pisier} G. Pisier. Martingales with values in uniformly
convex spaces. {\em Israel J. Math.} 20, no. 3--4, 326--350
(1975).


\bibitem{pisier-annals} G. Pisier. Holomorphic semigroups and the
geometry of Banach spaces. {\em Ann. Math. (2)}, 115, no. 2,
375--392 (1982).


\bibitem{pisier-varenna} G. Pisier. Probabilistic methods in the
geometry of Banach spaces. {\em Probability and Analysis (Varenna
1985)}, Springer Lect. Notes Math. 1206, 167--241 (1986).

\bibitem{pisier-book} G. Pisier. The Volume of Convex Bodies and
Banach Space Geometry. Cambridge University Press (1989).

\bibitem{pisier-xu} G. Pisier and Q. Xu. Random series in the real
interpolation spaces between the spaces $v_p$. {\em Geometrical
Aspects of Functional Analysis}, Springer Lect. Notes 1267,
185--209 (1986).

\bibitem{pitman} J. W. Pitman. One dimensional Brownian motion and
the three dimensional Bessel process. {\em Adv. Appl. Prob.} 7
511--526 (1975).

\bibitem{tsarkov}I. G. Tsar'kov. Extension of Hilbert-valued Lipschitz
mappings. {\em Vestik Moskov. Univ. Ser. I Math. Mekh.}, 6, 9--16
(1999).

\bibitem{vaisala} J. V\"ais\"al\"a. Gromov Hyperbolic Spaces.
Available at http://www.helsinki.fi/\string~jvaisala/grobok.pdf (2004).

\bibitem{ww} J. H. Wells and L. R. Williams. Embeddings and
Extensions in Analysis. Ergebnisse der Mathematik und ihrer
Grenzgebiete, Band 84, Springer-Verlag, New York (1975).

\bibitem{woj} P. Wojtaszczyk. Banach Spaces for Analysts.
Cambridge University Press, 1991.
\end{thebibliography}
\end{document}